\documentclass[12pt,a4paper]{article}
\usepackage{amsmath, amssymb, latexsym}
\usepackage{amsthm}
\newcommand{\N}{\mbox{$\mathbb N$}}
\newcommand{{\Z}}{{\bf Z}}
\newcommand{\R}{\mbox{$\mathbb R$}}

\newcommand{\omfp}{(\mbox{$\Omega$,$\cal F,$$P)$}}

\newtheorem{theorem}{Theorem}
\newtheorem{cor}[theorem]{Corollary}
\newtheorem{defn}[theorem]{Definition}
\newtheorem{lemma}[theorem]{Lemma}
\newtheorem{prop}[theorem]{Proposition}
\newtheorem{ex}[theorem]{Example}
\newtheorem{remark}[theorem]{Remark}

\newtheorem{assumption}[theorem]{Assumption}

\numberwithin{theorem}{section}
\numberwithin{equation}{section}

\newcommand{\beq}{\begin{equation}}
\newcommand{\eeq}{\end{equation}}
\newcommand{\beqo}{\begin{equation*}}
\newcommand{\eeqo}{\end{equation*}}
\newcommand{\bdm}{\begin{displaymath}}
\newcommand{\edm}{\end{displaymath}}
\newcommand{\beqar}{\begin{eqnarray}}
\newcommand{\eeqar}{\end{eqnarray}}
\newcommand{\beqaro}{\begin{eqnarray*}}
\newcommand{\eeqaro}{\end{eqnarray*}}
\newcommand{\bal}{\begin{align}}
\newcommand{\eal}{\end{align}}

\newcommand{\btheorem}{\begin{theorem}}
\newcommand{\etheorem}{\end{theorem}}
\newcommand{\blemma}{\begin{lemma}}
\newcommand{\elemma}{\end{lemma}}
\newcommand{\bremark}{\begin{remark}}
\newcommand{\eremark}{\end{remark}}
\newcommand{\bcor}{\begin{cor}}
\newcommand{\ecor}{\end{cor}}
\newcommand{\bex}{\begin{ex}}
\newcommand{\eex}{\end{ex}}
\newcommand{\bdefn}{\begin{defn}}
\newcommand{\edefn}{\end{defn}}
\newcommand{\bprop}{\begin{prop}}
\newcommand{\eprop}{\end{prop}}
\newcommand{\bproof}{\begin{proof}}
\newcommand{\eproof}{\end{proof}}
\newcommand{\bass}{\begin{assumption}}
\newcommand{\eass}{\end{assumption}}

\newcommand{\M}{\text{{\sc M}}}

\newcommand{\X}{{\bf X}}
\newcommand{\Y}{{\bf Y}}
\begin{document}

\title{Alternating birth-death processes
}
\author{Hans Daduna 
\thanks{daduna@math.uni-hamburg.de, ORCID 0000-0001-6570-3012, Daduna-AlternBDPr-20200511}
\\
    {  Hamburg University,
         Department of Mathematics},\\
          Bundesstrasse 55,
          20146 Hamburg,
          Germany
}
\maketitle

\begin{abstract}
We consider a continuous time Markov process on $\N_0$ which can be interpreted as generalized alternating birth-death process in a non-autonomous random environment. Depending on the status of the environment the process either increases until the environment changes and the process starts to decrease until the environment changes again, and the process restarts to increase, and so on, or its starts decreasing, reversing its direction due to environmental changes, et cetera.
The birth and death rates depend on the state (height, population size) of the birth-death process and the environment's transition rates depend on the state  of the birth-death process as well. Moreover, a birth or death event may trigger an immediate change of the environment.
Our main result is an explicit expression for the stationary distribution if the system is ergodic,
providing ergodicity conditions as well.

Removing the reflecting boundary at zero we obtain a two-sided version on $\mathbb{Z}$ of this alternating birth-death process, which for suitable parameter constellations is ergodic as well. We determine the stationary distribution.
This two-sided version is a  locally inhomogeneous discrete space version of the classical telegraph process.

We demonstrate that alternating birth-death processes in a random environment provide a versatile class of models from different areas of applications. Examples from the literature are discussed.
\end{abstract}
{\it MSC2020-classification}: 60J27 (Primary), 60J28, 60K25, 60K37 (Secondary)
\section{Introduction}\label{sect:intro}
Birth-death processes are used to describe the evolution of  different systems under random influences, e.g. population dynamics, queueing systems in various application areas, inventories.
In many situations occurring in these areas the systems under observation are influenced by an external environment which is subject to random influences as well. There exist a rich literature on birth-death processes in a random environment, where typically the environment develops autonomously as a Markov process of its own, i.e. the environment is autonomous. Queueing systems as special birth-death processes in an autonomous random environment are investigated e.g. by Zhu \cite{zhu:94}, Economou \cite{economou:05}, Tsitsiashvili, Osipova, Koliev, Baum \cite{tsitsiashvili;osipova;koliev;baum:02}, and Balsamo, Marin
\cite{balsamo;marin:13}.
These systems are typical examples of ``quasi-birth-death processes'' (QBD's) where the ``level'' represents the population size while the associated ``phase'' represents the environment. The level-independent QBDs model birth-death processes in an autonomous environment while the level-dependent QBDs are related to systems under  perturbations by non-autonomous environments, for classical references
see Chapter 6 (Queues in a Random Environment) in \cite{neuts:81}, and Example C in \cite{neuts:89}[p. 202, 203].\\

Our focus is on modified birth-death processes in a non-autonomous random environment. This is motivated by
several classes of applied probability models which are investigated in the literature, where
the environment is subject to perturbations originating from the dynamics of the birth-death process.
A class of examples represents models for integrated production-inventory systems, where the production facility is modeled as a queue, while the inventory represents the environment of that queue, a review is
\cite{krishnamoorthy;lakshmy;manikandan:11}, more recent results are a in \cite{krenzler;daduna:15}.
A different area where typically the environment of a queue is non-autonomous is the fields of wireless sensor networks, where the principle of a ``referenced node'' is used to aggregate the network and to consider only a single node, resp. its message queue and to incorporate the external conditions and the other nodes of the network into its environment. Sending and receiving messages by the referenced message queue (the ``birth-death process'') changes the status of the environment, for a discussion and an elaborated example see \cite{krenzler;daduna:14}.\\

While in standard birth-death processes in a random environment in any state (population size) $> 0$ births \underline{and} deaths occur with intensities which may depend on state and environment,
the alternating birth-death processes in a random environment 
considered here are allowed to move up (birth) or down (death) only when the environment is in the respective state, called $b$ or $d$.
Our interest in these class of processes originates from different models in the literature which
can be subsumed under a common scheme. Typical examples:\\
Kella and Whitt \cite{kella;whitt:92} investigated the structure of an infinite-capacity storage model
(where the maximal service capacity is greater than the arrival intensity) influenced by a two-state random environment. The environment alternates between states ``up'' and ``down''.
To investigate this model is motivated by queueing systems with random  interruptions  of service,
the status of the environment indicates the server's availability: If the environment
is down, no service is provided, only arrivals occur and consequently the content
of the storage increases pathwise according to some general stochastic process. If the environment is up, service is provided and arrivals occur and consequently the content decreases pathwise according to some general stochastic process because service intensity exceeds arrival intensity.\\
A rather general storage process which is closely related to the model in  \cite{kella;whitt:92}
is investigated by Boxma and Kella \cite{boxma;kella:14}. They consider a process which is either in up or in down state. In up-state it behaves like a L\'{e}vy process with no negative jumps and negative drift, while in down-state it behaves like a subordinator.

Consider now for the system in \cite{kella;whitt:92} the situation where during up times the content increases according to a Poisson process
 with intensity $\lambda$, while during down times the content decreases according to a
 Poisson process with intensity $\mu - \lambda > 0$ (notation of \cite{kella;whitt:92}).
Then the content process may be described  via discrete state space and can be considered as an ``alternating birth-death
process in a random environment'' which increases by births as long as the environment is down, and decreases by deaths as long as the environment is up.

Our focus in this paper is in a first step on  similar processes on state space $\N_0$ with general transition mechanism: 
Birth and death rates depend on the state
(height, population size) of the birth-death process (and the status of the environment), and the environment's transition rates depend on the state  of the birth-death process. Furthermore, an immediate change of the environment status may be triggered by a birth or death occurring,
i.e. an arrival or departure in terms of queueing models.

In a second step we remove the reflecting boundary at zero and extend the process to a two-sided birth-death process on state space $\mathbb Z$.
This results in a process which makes jumps on $\mathbb Z$ either to the left for a random time duration  and then suddenly changes the direction of his moves to the right, until the next change.
The duration of times with jumps of constant direction is in our setting controlled by the status of the environment. On state space $\R$
such a process is known as generalized telegraph process, for a short introduction see
the fundamental paper of Kac \cite{kac:74} and the more recent work of Stadje and Zacks \cite{stadje;zacks:04}, de Gregorio \cite{degregorio:10}, and Crimaldi, di Crescenzo, Iuliano, and Martinucci 
\cite{crimaldi;dicrescenzo;iuliano;martinucci:13}. The simplest telegraph processes describe the motion of a particle on the continuous real line with constant speed, the direction is reversed at jump times of a
Poisson process. We remark that Kac derives the equation which governs the density of the particle's position starting from a discrete version on the one dimensional lattice $\mathbb Z$ of the particle's walk \cite{kac:74}.\\
A generalization of the telegraph process is defined by Ratanov \cite{ratanov:20}: The particle moves irregularly governed by different L\'evy processes and the switching between these processes is governed by an underlying two-state Markov process.
(In our setting the L\'evy processes are Poisson.)
More related processes which generalize the original telegraph processes can be found in the references there.

Under the heading of ``Double-ended queue'' and ``unrestricted linear random walk'' such two-sided birth-death processes have found many applications.
An early survey of the two-sided $M/M/1/\infty$ queue denoted by $\infty^2/M/M$, resp. two-sided birth-death process on state space $\mathbb Z$ is of Conolly \cite{conolloy:71}. More details can be found in the book of Srivastava and Kashyap \cite{srivastava;kashyap:82}.
A recent in depth study is the paper of Liu, Gong, and Kulkarni \cite{liu;gong;kulkarni:15}.

From a general point of view our model with state space $\mathbb N_0\times\{b,d\}$ is related to the processes which Falin and Gomez-Corral introduced in \cite{falin;gomez-corral:00} as a class of  bivariate Markov processes. The construction of these processes is motivated by models from teletraffic analysis
and the resulting class of models encompass many different retrial systems from the previous literature.

\textbf{Structure of the paper.} In Section \ref{sect:ModelOneSided} we consider the one-sided birth-death process on $\mathbb{N}_0$, compute the stationary distribution for the ergodic system, and discuss several examples. Because the process is not reversible we
provide the stationary distribution for the finite space system
direct in Section \ref{sect:OneSidedFinite}. In Section \ref{sect:REGULAROneSided} we discuss regularity issues.
In Section \ref{sect:ModelTwoSided} we consider the two-sided birth-death process on $\mathbb{Z}$, compute the stationary distribution for the ergodic system, and discuss the telegraph process as example. We show that our approach enables us to stabilize the moving particle's process.
In \ref{sect:REGULARTwoSided} we discuss regularity.\\  
{\bf Assumptions and conventions.} We assume throughout that all processes which occur are defined on a common probability space on $\omfp$ and have cadlag paths.
By construction all processes which will occur are conservative,  i.e. their
intensity matrix (generator matrix) has zero row  sums.

\section{One-sided alternating birth-death process} \label{sect:ModelOneSided}

\begin{figure}
\begin{picture}(20,4)
\thicklines
\multiput(1,0)(4,0){4}{\circle{1.1}}
\multiput(1,4)(4,0){4}{\circle{1.1}}
\multiput(1.7,4)(4,0){3}{\vector(1,0){2.6}}
\multiput(4.3,0)(4,0){3}{\vector(-1,0){2.6}}
\multiput(0.7,0.6)(4,0){4}{\vector(0,1){2.9}}
\multiput(1.3,3.4)(4,0){4}{\vector(0,-1){2.9}}
\multiput(1.7,3.5)(4,0){3}{\vector(1,-1){2.9}}
\multiput(4.3,0.5)(4,0){3}{\vector(-1,1){2.9}}
\put(0.9,3.9){\makebox[2mm]{n-1,b}}
\put(0.9,-0.1){\makebox[2mm]{n-1,d}}
\put(4.9,3.9){\makebox[2mm]{n,b}}
\put(4.9,-0.1){\makebox[2mm]{n,d}}
\put(8.9,3.9){\makebox[2mm]{n+1,b}}
\put(8.9,-0.1){\makebox[2mm]{n+1,d}}
\put(12.9,3.9){\makebox[2mm]{n+2,b}}
\put(12.9,-0.1){\makebox[2mm]{n+2,d}}
\put(0.1,2){\makebox[2mm]{$\beta_{n-1}$}}
\put(1.7,2){\makebox[2mm]{$\delta_{n-1}$}}
\put(4.2,2){\makebox[2mm]{$\beta_{n}$}}
\put(5.6,2){\makebox[2mm]{$\delta_{n}$}}
\put(8.1,2){\makebox[2mm]{$\beta_{n+1}$}}
\put(9.7,2){\makebox[2mm]{$\delta_{n+1}$}}
\put(12.1,2){\makebox[2mm]{$\beta_{n+2}$}}
\put(13.2,2){\makebox[2mm]{$\delta_{n+2}$}}
\put(3,4.1){\makebox[2mm]{$\lambda_{n-1}$}}
\put(3,-0.3){\makebox[2mm]{$\mu_{n}$}}
\put(7,4.1){\makebox[2mm]{$\lambda_{n}$}}
\put(7,-0.3){\makebox[2mm]{$\mu_{n+1}$}}
\put(11,4.1){\makebox[2mm]{$\lambda_{n+1}$}}
\put(11,-0.3){\makebox[2mm]{$\mu_{n+2}$}}
\put(2.5,3.2){\makebox[2mm]{$\kappa_{n-1}$}}
\put(6.5,3.2){\makebox[2mm]{$\kappa_{n}$}}
\put(10.5,3.2){\makebox[2mm]{$\kappa_{n+1}$}}
\put(3.8,0.5){\makebox[2mm]{$\nu_{n}$}}
\put(7.5,0.5){\makebox[2mm]{$\nu_{n+1}$}}
\put(11.5,0.5){\makebox[2mm]{$\nu_{n+2}$}}
\multiput(0,0)(-0.5,0){3}{\circle*{0.2}}
\multiput(0,4)(-0.5,0){3}{\circle*{0.2}}
\multiput(14,0)(0.5,0){3}{\circle*{0.2}}
\multiput(14,4)(0.5,0){3}{\circle*{0.2}}
\end{picture}

\caption{Transition graph for the one-sided process($n>1$)}\label{fig:1-sided-transitions1}
\end{figure}
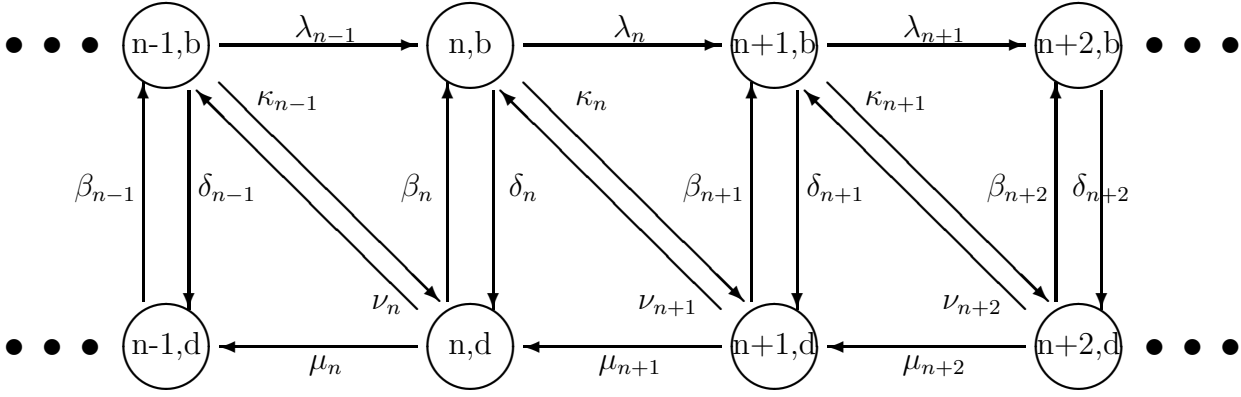

Let $\X=(X(t):t\geq 0)$ with $X(t): \omfp\to (\N_0,2^{\N_0})$
be a variant of the standard birth-death processes the evolution of which depends on an external environment which changes randomly between two states $\{b, d\}$. We denote the environment process by
$\Y=(Y(t):t\geq 0)$ with $Y(t): \omfp\to (\{b, d\},2^{\{b, d\}})$.
Whenever the environment is in state $b$ (indicating that {\em births} may occur) the process $\X$ moves upwards only, due to births occurring, and
whenever the environment is in state $d$ (indicating that {\em deaths} may occur) the process $\X$ moves downwards only, due to deaths occurring.
The joint process with state space $\N_0\times \{b, d\}$ is denoted by $\Z=(\X,\Y)$ with
$Z(t) = (X(t),(Y(t))$. Movements of $\Z$ are governed by transition intensity matrix
 $Q=(q(z,z'):z,z'\in \N_0\times \{b, d\})$ with strictly positive entries as follows for $n\geq 0$:
\beqaro
q(n,b;n+1,b) &=& \lambda_n,~~\text{a birth occurs when the population size is}~~n, \\ 
q(n,d;n-1,d) &=& \mu_n 1_{(n>0)},~~\text{a death occurs when the population size  is}~~n,\\ 
q(n,b;n,d) &=& \delta_n,~~\text{environment changes from ``birth'' to ``death'' when}\nonumber\\ &&\text{the population size is}~~n, \\ 
q(n,d;n,b) &=& \beta_n,~~\text{environment changes from ``death'' to ``birth'' when}\nonumber\\
 &&\text{the population  size is}~~n,\\ 
q(n,b;n+1,d) &=& \kappa_n,~~\text{a birth occurs {\bf and} environment changes from ``birth''}\nonumber\\
 &&\text{ to ``death'' when the population size is}~~n, \\ 
q(n,d;n-1,b) &=& \nu_n 1_{(n>0)},~~\text{a death occurs {\bf and} environment changes from }\nonumber\\
 &&\text{``death'' to ``birth'' when the population  size is}~~n. 
\eeqaro
The diagonal elements of $Q$ are chosen such that row sums are zero. Unless otherwise indicated for special situations we assume throughout that the parameters
$\lambda_n, \kappa_n, \delta_n, \mu_n, \nu_n, \beta_n$ are strictly positive. A typical scenario of the transition graph is visualized in Figure \ref{fig:1-sided-transitions1}.

 \unitlength1cm

\subsection{Stationary distribution for the one-sided process}\label{sect:GBOneSided}
We assume in this section that the birth-death process is non-exploding, i.e. is regular. (In Section \ref{sect:REGULAROneSided} we will discuss this in more detail.)
The global balance equation for $\Z$ are with unknown values $x(n,a)$ for
$(n,a)\in \mathbb N_0\times\{b,d\}$
\beqar
x(0,b)(\lambda_0 + \kappa_0 + \delta_0) &=& x(1,d) \nu_1 + x(0,d) \beta_0 ,\label{eq:GBOne1}\\
x(0,d) \beta_0 &=& x(1,d) \mu_1+ x(0,b) \delta_0, \label{eq:GBOne2}\\
\text{and for }~~n\geq 1&&\nonumber\\
x(n,b)(\lambda_n + \kappa_n + \delta_n) &=& x(n-1,b)\lambda_{n-1} + x(n+1,d) \nu_{n+1}
+ x(n,d) \beta_n, \label{eq:GBOne3}\\
x(n,d)(\mu_n + \nu_n + \beta_n) &=& x(n+1,d)\mu_{n+1} +  x(n-1,b) \kappa_{n-1} +
x(n,b) \delta_{n}. \label{eq:GBOne4}
\eeqar

\btheorem\label{thm:OneSidedStationary}
Denote for $n\geq 0$ by $\Lambda_n := \lambda_n + \kappa_n$ the total uprate out of $n$ and $\M_n := \mu_n + \nu_n$ the total downrate out of $n$.
The global balance equation for $\Z$ are solved with any $x(0,d)>0$ by
\beq
x(n,b)= x(0,d)\beta_0 \M_{n+1}
\frac{\prod_{k = 1}^{n}\left(\Lambda_{k-1} \beta_{k} +\M_{k}\lambda_{k-1}\right)}
{\prod_{k = 0}^{n}\left(\Lambda_{k} \mu_{k+1}+\M_{k+1}\delta_{k}\right)}, \quad n\geq 0,
\label{eq:Oneuptheo1}
\eeq
and
\beq
x(n,d) = x(0,d) \beta_0 \Lambda_{n-1}
\frac{\prod_{k = 1}^{n-1}\left(\Lambda_{k-1} \beta_{k} +\M_{k}\lambda_{k-1}\right)}
{\prod_{k = 0}^{n-1}\left(\Lambda_{k} \mu_{k+1}+\M_{k+1}\delta_{k}\right)}, \quad n\geq 1.
\label{eq:Oneuptheo2}
\eeq
$\Z$ is ergodic if and only if
\beqo 
C(1) = \sum_{n=0}^\infty
\frac{\prod_{k = 1}^{n}\left(\Lambda_{k-1} \beta_{k} +\M_{k}\lambda_{k-1}\right)}
{\prod_{k = 0}^{n}\left(\Lambda_{k} \mu_{k+1}+\M_{k+1}\delta_{k}\right)}
\cdot(\Lambda_n+\M_{n+1}) < \infty
\eeqo
\etheorem
\bproof
The key observation is that a sequence of simple cuts reduces the computation to investigating a linear structure. We have for all  $n\geq 0$ the cut equation
\beqar
x(n,b)(\lambda_n + \kappa_n) &=&  x(n+1,d)(\mu_{n+1} + \nu_{n+1}), \nonumber\\ 
\text{which yields for all}~~ n\geq 0 ~~~~~~~
x(n+1,d) &=& x(n,b)\frac{\Lambda_n}{\M_{n+1}} .  \label{eq:cutOne2}
\eeqar
This allows to substitute in a first step the variables $x(n,d)$. In \eqref{eq:GBOne2} we obtain
\beqar
x(0,d) \beta_0 &=& [x(0,b)\frac{\Lambda_0}{\M_{1}}] \mu_1+ x(0,b) \delta_0,\nonumber\\
\text{and therefore}~~~x(0,b) &=& x(0,d)\beta_0 \left(\frac{\Lambda_0}{\M_{1}}\mu_1+  \delta_0\right)^{-1}.
\label{eq:substOne1}
\eeqar
From \eqref{eq:GBOne3} we obtain for $n\geq 1$
\beqaro
x(n,b)(\lambda_n + \kappa_n + \delta_n) = x(n-1,b)\lambda_{n-1} +
[x(n,b)\frac{\Lambda_n}{\M_{n+1}}] \nu_{n+1}
+ [x(n-1,b)\frac{\Lambda_{n-1}}{\M_{n}}]  \beta_n,&&\\
\text{and so}~~~
x(n,b)\frac{(\lambda_n + \kappa_n + \delta_n)\M_{n+1} - \Lambda_n \nu_{n+1}}{\M_{n+1}}
= x(n-1,b)\frac{(\lambda_{n-1}\M_n + \Lambda_{n-1} \beta_n)}{\M_{n}},&&
\eeqaro
which leads to
\beqo
x(n,b)= x(n-1,b)\frac{\M_{n+1}}{\M_{n}}\cdot
\frac{\Lambda_{n-1} \beta_{n} +\M_{n}\lambda_{n-1}}{\Lambda_{n} \mu_{n+1}+\M_{n+1}\delta_{n}},\quad n\geq 1.
\eeqo
Iterating and using \eqref{eq:substOne1} we obtain for all  $n\geq 0$
\beqar
x(n,b)&=& x(0,d)\beta_0 \left(\frac{\Lambda_0}{\M_{1}}\mu_1+  \delta_0\right)^{-1}
\cdot \prod_{k = 1}^{n}\frac{\M_{k+1}}{\M_{k}}\cdot
\prod_{k = 1}^{n}
\frac{\Lambda_{k-1} \beta_{k} +\M_{k}\lambda_{k-1}}{\Lambda_{k} \mu_{k+1}+\M_{k+1}\delta_{k}}\nonumber\\
&=& x(0,d)\beta_0 \left(\frac{\Lambda_0}{\M_{1}}\mu_1+  \delta_0\right)^{-1}
\cdot \frac{\M_{n+1}}{\M_{1}}\cdot
\prod_{k = 1}^{n}
\frac{\Lambda_{k-1} \beta_{k} +\M_{k}\lambda_{k-1}}{\Lambda_{k} \mu_{k+1}+\M_{k+1}\delta_{k}}\nonumber\\
&=& x(0,d)\beta_0 \M_{n+1}
\frac{\prod_{k = 1}^{n}\left(\Lambda_{k-1} \beta_{k} +\M_{k}\lambda_{k-1}\right)}
{\prod_{k = 0}^{n}\left(\Lambda_{k} \mu_{k+1}+\M_{k+1}\delta_{k}\right)}.
\label{eq:upOne1}
\eeqar
Adding \eqref{eq:GBOne1}  and \eqref{eq:GBOne2} yields
\beqo
x(1,d)= x(0,b) \frac{\Lambda_0}{\M_1},
\eeqo
and applying \eqref{eq:substOne1} we obtain
\beq\label{eq:downOne1}
x(1,d)= x(0,d) \beta_0 \Lambda_0 \frac{1}{\Lambda_0\mu_1+\M_1 \delta_0},
\eeq
and for $n\geq 1$ we obtain from \eqref{eq:cutOne2} and \eqref{eq:upOne1}
\beqo
x(n,d) = x(0,d) \beta_0 \Lambda_{n-1}
\frac{\prod_{k = 1}^{n-1}\left(\Lambda_{k-1} \beta_{k} +\M_{k}\lambda_{k-1}\right)}
{\prod_{k = 0}^{n-1}\left(\Lambda_{k} \mu_{k+1}+\M_{k+1}\delta_{k}\right)},
\eeqo
which indeed is in line with \eqref{eq:downOne1}.
\eproof
\begin{remark}
\textbf{(a)} The ergodicity criterion $C(1)<\infty$ is well-suited for detailed investigation via ``Extensions of the Bertrand-De Morgan test'' for convergence of series with positive summands. For an indepth study see \cite{abramov:20}.\\
\textbf{(b)}	
A stationary  birth-death process is reversible for any parameter setting. This is not the case for the alternating birth-death process in a random environment, which can be seen easily by writing down the transition intensities of its time reversal. 
\end{remark}

\subsubsection{Example: $M/M/1/\infty$ retrial queues with general retrial policy}\label{sect:RetrialQueue}
Retrial queues are classical models for telephony, our description follows \cite{artalejo;gomez-corral:08}[Section 2]:
Customers arrive in a Poisson$-\lambda$ stream at a single server without waiting room. Service times are
$\exp(\delta)$-distributed. The server's states are {\em empty} $= 0$ or {\em busy} $= 1$.\\
If an arriving customer finds the server empty he enters the system and  service starts immediately.
If an arriving customer finds the server busy he leaves the service area and enters the so-called orbit
where he successively retries to enter the server. If on his retrial he finds the server idling
he enters the system and  service starts immediately. Otherwise he stays-on in the orbit.
There are various retrial policies described in the literature, characterized mainly by the retrial intensities, see \cite{falin:13}[p. 417].

\begin{figure}
	\begin{picture}(20,4)
	\thicklines
	\multiput(1,0)(4,0){4}{\circle{1.1}}
	\multiput(1,4)(4,0){4}{\circle{1.1}}
	\multiput(1.7,4)(4,0){3}{\vector(1,0){2.6}}
	\multiput(0.7,0.6)(4,0){4}{\vector(0,1){2.9}}
	\multiput(1.3,3.4)(4,0){4}{\vector(0,-1){2.9}}
	\multiput(4.3,0.5)(4,0){3}{\vector(-1,1){2.9}}
	\put(0.9,3.9){\makebox[2mm]{n-1,1}}
	\put(0.9,-0.1){\makebox[2mm]{n-1,0}}
	\put(4.9,3.9){\makebox[2mm]{n,1}}
	\put(4.9,-0.1){\makebox[2mm]{n,0}}
	\put(8.9,3.9){\makebox[2mm]{n+1,1}}
	\put(8.9,-0.1){\makebox[2mm]{n+1,0}}
	\put(12.9,3.9){\makebox[2mm]{n+2,1}}
	\put(12.9,-0.1){\makebox[2mm]{n+2,0}}
	\put(0.1,2){\makebox[2mm]{$\lambda$}}
	\put(1.7,2){\makebox[2mm]{$\delta$}}
	\put(4.2,2){\makebox[2mm]{$\lambda$}}
	\put(5.6,2){\makebox[2mm]{$\delta$}}
	\put(8.1,2){\makebox[2mm]{$\lambda$}}
	\put(9.7,2){\makebox[2mm]{$\delta$}}
	\put(12.1,2){\makebox[2mm]{$\lambda$}}
	\put(13.4,2){\makebox[2mm]{$\delta$}}
	\put(3,4.1){\makebox[2mm]{$\lambda$}}
	\put(7,4.1){\makebox[2mm]{$\lambda$}}
	\put(11,4.1){\makebox[2mm]{$\lambda$}}
	\put(2.7,1){\makebox[2mm]{$\nu_{n}$}}
	\put(6.4,1){\makebox[2mm]{$\nu_{n+1}$}}
	\put(10.4,1){\makebox[2mm]{$\nu_{n+2}$}}
	\multiput(0,0)(-0.5,0){3}{\circle*{0.2}}
	\multiput(0,4)(-0.5,0){3}{\circle*{0.2}}
	\multiput(14,0)(0.5,0){3}{\circle*{0.2}}
	\multiput(14,4)(0.5,0){3}{\circle*{0.2}}
	\end{picture}
	
	\caption{Transition graph of the retrial queue with general retrial policy
		\cite{falin:13}}\label{fig:RetrialQueue2}
	
\end{figure}
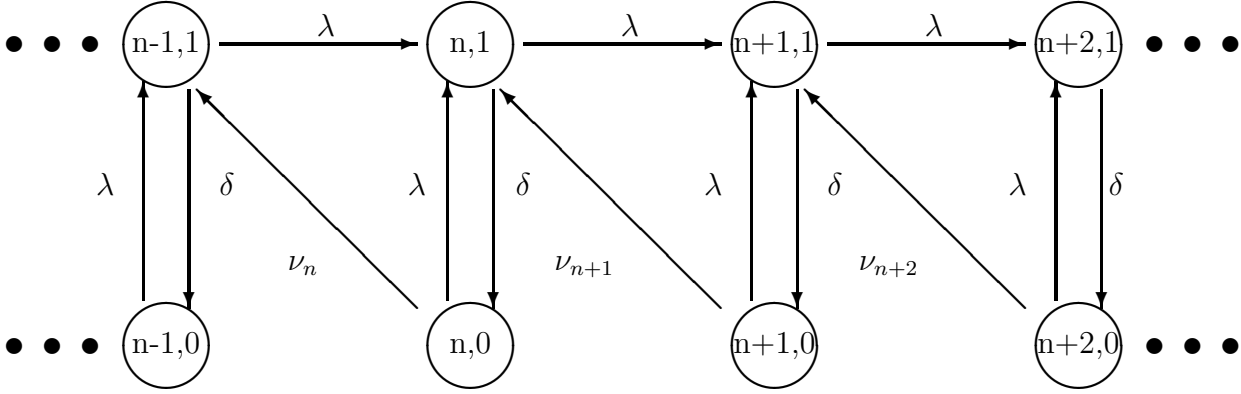

A general retrial policy is as follows: Customers in the orbit queue up according to a
First-Come-First-Served (FCFS) regime.
If there are $n\geq 1$ customers in the orbit the overall retrial intensity is $\nu_{n}$, i.e, if at time $t\geq 0$ the orbit population size is $n\geq 1$ then the one customer at the head of the queue attempts to enter the server
during $[t, t+\Delta)$ with probability $\nu_n\cdot \Delta + o(\Delta)$. If the server is free ( $= 0$)
it becomes occupied and the orbit population size changes to $n-1$. If the server is occupied  ($= 1$)
the attempt is not successful and the orbit population size stays at $n$, with customers in the same order, the transition diagram of Falin's retrial queue with general retrial scheme is shown in Figure \ref{fig:RetrialQueue2}.

\begin{figure}
	\unitlength1cm
	\begin{picture}(20,4.3)
	\put(0,1){\vector(1,0){2.8}}
	\put(0,0.6){Arrival stream}
	\put(2.9,1){\circle{0.2}}
	\put(3,1){\vector(1,0){2.8}}
	\put(3.2,1.2){admitted if}
	\put(3.2,0.6){server idle}
	\put(2.9,1.1){\vector(0,1){2}}
	\put(8,1){\vector(1,0){2.5}}
	\put(8.2,0.6){Departure stream}
	\put(6,0.5){\framebox(2,1){{\circle*{0.2}}}}
	\put(6.3,0.1){\mbox{Server}}
	\put(6.3,-0.3){\mbox{occupied}}
	\put(4,3.7){\oval(5,1)}
	\put(5.5,3){\vector(0,-1){1.8}}
	\put(2,4.4){ORBIT}
	\put(4.2,2.5){Retrial}
	\put(4.2,2.1){stream}
	\put(1.4,2.5){Blocked}
	\put(1.1,2){customers}
	\multiput(4,3.7)(0.4,0){5}{\circle*{0.2}}
	\multiput(0.7,1)(0.5,0){4}{\circle*{0.2}}
	\multiput(8.7,1)(0.5,0){3}{\circle*{0.2}}
	\end{picture}\\
	\caption{The retrial queue}\label{fig:RetrialQueue1}
\end{figure}
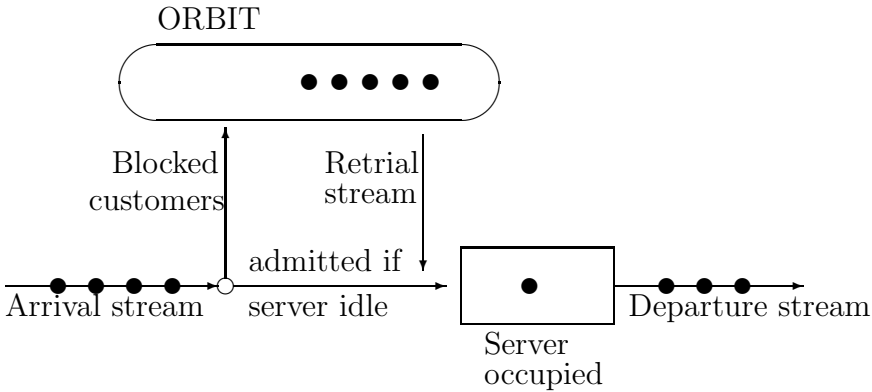

This system fits into the model of Section \ref{sect:ModelOneSided} by interpreting the server as the
orbit's environment with $b\leftrightarrow 1$ and $d\leftrightarrow  0$.
The orbit's population size is the state of the
birth-death process which (only!) increases with rate $\lambda$ as long as the environment's state is $1$ (new arrivals are send to the orbit because the server is occupied).
The orbit's population size (only!) decreases with rate $\nu_n$ as long as the environment's state is $0$
and the ``queue length'' of the orbit is $n>0$. If this is the case, a downward jump of the orbit's queue length
from $n$ to $n-1$ triggers an immediate change of the environment from $0$ to $1$: The server is occupied now and the orbit size starts to increase again.

According to \cite{falin:13}[p. 417] of special interest are the following retrial policies:\\
$\circ$  constant  policy, i.e. $\nu_n = \alpha, \forall n\geq 1$; the orbit is organized as a standard $M/M/1/\infty$ queue under First-Come-First-Served (FCFS);\\
$\circ$  classical policy, i.e.  $\nu_n = \nu\cdot n, \forall n\geq 1$; the orbit is organized as a standard $M/M/\infty$ queue;\\
$\circ$  linear  policy, i.e.  $\nu_n = \alpha + \nu\cdot n, \forall n\geq 1$.\\
Details and more references can be found in \cite{falin:13}
and \cite{artalejo;gomez-corral:08}.

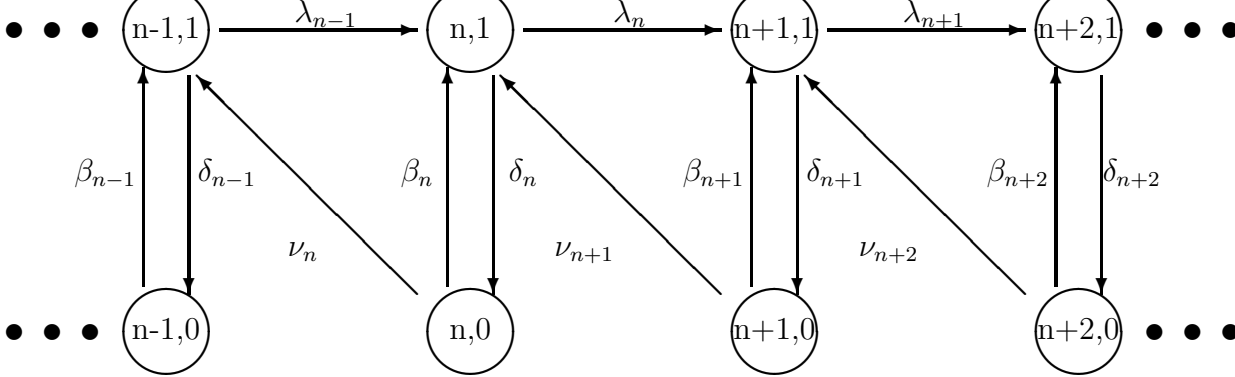
\begin{figure}
	\begin{picture}(20,4)
	\thicklines
	\multiput(1,0)(4,0){4}{\circle{1.1}}
	\multiput(1,4)(4,0){4}{\circle{1.1}}
	\multiput(1.7,4)(4,0){3}{\vector(1,0){2.6}}
	\multiput(0.7,0.6)(4,0){4}{\vector(0,1){2.9}}
	\multiput(1.3,3.4)(4,0){4}{\vector(0,-1){2.9}}
	\multiput(4.3,0.5)(4,0){3}{\vector(-1,1){2.9}}
	\put(0.9,3.9){\makebox[2mm]{n-1,1}}
	\put(0.9,-0.1){\makebox[2mm]{n-1,0}}
	\put(4.9,3.9){\makebox[2mm]{n,1}}
	\put(4.9,-0.1){\makebox[2mm]{n,0}}
	\put(8.9,3.9){\makebox[2mm]{n+1,1}}
	\put(8.9,-0.1){\makebox[2mm]{n+1,0}}
	\put(12.9,3.9){\makebox[2mm]{n+2,1}}
	\put(12.9,-0.1){\makebox[2mm]{n+2,0}}
	\put(0.1,2){\makebox[2mm]{$\beta_{n-1}$}}
	\put(1.7,2){\makebox[2mm]{$\delta_{n-1}$}}
	\put(4.2,2){\makebox[2mm]{$\beta_{n}$}}
	\put(5.6,2){\makebox[2mm]{$\delta_{n}$}}
	\put(8.1,2){\makebox[2mm]{$\beta_{n+1}$}}
	\put(9.7,2){\makebox[2mm]{$\delta_{n+1}$}}
	\put(12.1,2){\makebox[2mm]{$\beta_{n+2}$}}
	\put(13.6,2){\makebox[2mm]{$\delta_{n+2}$}}
	\put(3,4.1){\makebox[2mm]{$\lambda_{n-1}$}}
	\put(7,4.1){\makebox[2mm]{$\lambda_{n}$}}
	\put(11,4.1){\makebox[2mm]{$\lambda_{n+1}$}}
	\put(2.7,1){\makebox[2mm]{$\nu_{n}$}}
	\put(6.4,1){\makebox[2mm]{$\nu_{n+1}$}}
	\put(10.4,1){\makebox[2mm]{$\nu_{n+2}$}}
	\multiput(0,0)(-0.5,0){3}{\circle*{0.2}}
	\multiput(0,4)(-0.5,0){3}{\circle*{0.2}}
	\multiput(14,0)(0.5,0){3}{\circle*{0.2}}
	\multiput(14,4)(0.5,0){3}{\circle*{0.2}}
	\end{picture}
	
	\caption{Transition graph of the retrial queue with general retrial policy and state dependent arrival, service, and 
		retrial rates \cite{falin;gomez-corral:00}} \label{fig:RetrialQueue3}
	
\end{figure}
Falin and Gomez-Corral  \cite{falin;gomez-corral:00} considered the retrial queue from Figure \ref{fig:RetrialQueue2} in a more versatile version.
The arrival rates are different for different size of the orbit queue and depend furthermore on whether the server is busy or idling. The service rate depends on the size of the orbit queue. 

This results in a general ``bivariate Markov process arising in the theory of single-server retrial queues'' in \cite{falin;gomez-corral:00}, which is closer to our general system.
With the notation of  Section \ref{sect:ModelOneSided} the joint ``(orbit population size/server status)'' process  $\Z=(\X,\Y)$ is governed by the following intensity matrix.
\beqaro
q(n,1;n+1,1) &=&  \lambda_n~~\text{when the server is busy and the orbit population size is}~n,\nonumber\\
&&  \text{customers arrive with intensity}~
\lambda_n~~\text{and enter the orbit  }\\ 
q(n,0;n,1) &=& \beta_n,~~\text{when the server is idling and the orbit population size is}~n,\nonumber\\
&&  \text{customers arrive with intensity}~
\beta_n~~\text{and enter the orbit}\\ 
q(n,1;n,0) &=& \delta_n~~~~\text{when the server is busy and the orbit population size is}~n,\nonumber\\
&&\text{ service is provided with rate}  \delta_n
\\ 
q(n,0;n-1,1) &=& \nu_n 1_{(n>0)},~~\text{when the server is idling and the orbit population size is}~n,\nonumber\\
&&\text{the customer at the head of the orbit queue retries successfully}\nonumber\\
&&\text{ with rate}~\nu_{n}~\text{leaves the orbit   and enters the server}
\eeqaro
The structure of the transition intensity graph is presented in Figure  \ref{fig:RetrialQueue3}.
Assuming that all the rates depicted in Figure \ref{fig:RetrialQueue3} are positive, the associated queue lengths processes for the server-orbit system is irreducible. In case of positive recurrence the stationary distribution is with suitably scaled $x(0,d)$
\beq\label{eq:RetrialGeneral-b}
x(n,b) = x(0,d)\frac{\beta_0}{\delta_0}
\left(\prod_{k=1}^{n}\frac{\lambda_{k-1}}{\delta_{k}}\right)
\left(\prod_{k=1}^{n}\frac{\beta_{k}+\nu_{k}}{\nu_{k}}\right),
\eeq
and
\beq\label{eq:RetrialGeneral-d}
x(n,d) = x(0,d)\frac{\beta_0}{\delta_0}
\left(\prod_{k=1}^{n-1}\frac{\lambda_{k-1}}{\delta_{k}}\right)
\left(\prod_{k=1}^{n-1}\frac{\beta_{k}+\nu_{k}}{\nu_{k}}\right)
\frac{\lambda_{n-1}}{\nu_{n}}.
\eeq
For the case of general retrial policy with state independent arrival and service rates from Figure \ref{fig:RetrialQueue2} this boils down to Falin's result \cite{falin:13}[formulas (3), (2)].
\beq\label{eq:RetrialFalin-b}
x(n,b) = x(0,d) \left(\frac{\lambda}{\delta}\right)^{n+1}
\left(\frac{(\lambda+\nu_{1})\cdots(\lambda+\nu_{n})}
{\nu_{1}\cdots\nu_{n}}\right),
\eeq
and (with $\nu_{0}=0$)
\beq\label{eq:RetrialFalin-d}
x(n,d) = x(0,d)  \left(\frac{\lambda}{\delta}\right)^{n}
\left(\frac{(\lambda+\nu_{0})(\lambda+\nu_{1})\cdots(\lambda+\nu_{n})}{\nu_{1}\cdots\nu_{n}}\right).
\eeq
The formulas \eqref{eq:RetrialGeneral-b} for \eqref{eq:RetrialFalin-b}, resp. \eqref{eq:RetrialGeneral-d} for \eqref{eq:RetrialFalin-d}, reveal interesting structural properties of the formulas (3), (2) of Falin via expressions with state dependent rates: E.g. the factor $\left(\frac{\lambda}{\delta}\right)^{n+1}$ is composed from rates of different character, as can be seen from
$\frac{\beta_0}{\delta_0}
\left(\prod_{k=1}^{n}\frac{\lambda_{k-1}}{\delta_{k}}\right)$.
Similar insights are given for the other terms.

\subsubsection{Example: A dam model}\label{sect:dam}

Dam models have been dealt with as special storage models in the literature. A simple model of a dam is described by Prabhu \cite{prabhu:80}[p.71] as a model for a water reservoir with effectively infinite capacity. The input of the dam is a L\'{e}vy process
while the release from the dam is at unit rate except when the dam is empty.
The release of water is controlled by opening, resp. closing a gate.
In case of constant inflow this results in a dam model with constant input flow and gated outflow:\\
$\bullet$ Input flow  is with rate $\lambda$.\\
$\bullet$ Maximal  outflow rate when the gate is opened is $\theta$.\\
$\bullet$ So, controlled rate of decrease is $\theta - \lambda$ if the controller opens the gate (outflow - inflow).

If the time duration of open gate is exponential-$\beta$ and the closed-gate times are
exponential-$\delta$ we have a continuous state alternating fill-release process as described in Section
\ref{sect:ModelOneSided} with homogeneous-in-space transition structure.
The controller in this setting can be considered as the environment of the dam.\\
The connection to the alternating birth-death process: If the input and the output are discretized as Poisonian flows we can apply the above results e.g for determining equilibrium behaviour.
It turns out that the system is stable if and only if the natural condition 
$
\lambda/\delta < (\theta - \lambda)/\beta
$
is fulfilled. The stationary distribution $\pi$ of the stable system is  
\beqaro
\pi(0,d)&=&\frac{(\theta-\lambda)-\lambda \beta}{(\theta-\lambda)(\beta+\delta)},\\
\pi(n,d)&=& \pi(0,d) \frac{\beta\lambda}{(\lambda + \delta)(\theta-\lambda)}
\left(\frac{\lambda(\beta + \theta -\lambda)}{(\theta-\lambda)(\lambda+\delta)} \right)^{n-1},\quad n\geq 1,\\
\pi(n,b)&=& \pi(0,d) \frac{\beta}{\lambda + \delta}
\left(\frac{\lambda(\beta + \theta -\lambda)}{(\theta-\lambda)(\lambda+\delta)} \right)^n,\quad n\geq 0.
\eeqaro

\subsubsection{Example: Fluid queues}\label{sect:FluidQueue}
Adan and Resing \cite{adan;resing:96}[Section 1] investigated a fluid queue which resembles the dam model from the Section \ref{sect:dam}, but is in its simple version without control feature. The fluid queue with state space $[0,\infty)$ changes the queue size continuously: Leak rate is constant $=1$, while the input rate is determined by an alternating renewal process $(X_1,Y_1,X_2,Y_2,\dots)$:
During $X_i$-periods the input rate is $=2$, while during $Y_i$-periods the input rate is $=0$. So during $X_i$-periods the buffer content increases deterministically with rate $=1$ and decreases
deterministically during $Y_i$-periods with rate $=1$. It follows from the homogeneity assumptions of the system that the buffer content $Z_i$ at the beginning of the $i$th period $X_i$ follows the
waiting time recursion for the standard $G/G/1/\infty$ queue
\beqo
Z_{i+1} =\max(Z_i+X_i - Y_i,0).
\eeqo
With $X_i\sim \exp(\beta)$ and $Y_i\sim \exp(\delta)$ we have the dam model from the Section \ref{sect:dam}. Discretizing the buffer content and the inflow and out flow streams as Poisonian flows
connects this model with the alternating birth-death process 
formalism from Section \ref{sect:ModelOneSided}.\\
In \cite{adan;resing:96}[Section 2] the authors specialize the distributions of the alternating renewal process to be generated by the alternating busy $(X_i)$ and idle $(Y_i)$ periods of an $M/M/1/\infty$ queue. This can be considered as a system where an autonomous $M/M/1/\infty$ queue controls the opening of the dam.

In the follow-up paper \cite{adan;vandoorn;resing;scheinhardt:98},
Adan, van Doorn, Resing, and Scheinhardt
investigate a similar model under the assumption that  the fluid  queue and a controlling (generalized)  $M/M/1/N$ queue 
interact in both directions, because the environment-$M/M/1/N$
queue has service rates depending on the size of the fluid  queue, i.e. the environment-$M/M/1/N$ queue is non-autonomous.
This property is inherent in the very definition of the alternating birth-death process in a  random environment of Section \ref{sect:ModelOneSided} and \ref{sect:ModelTwoSided}.
For discussion of autonomous versus non-autonomous environment for the fluid queue see \cite{adan;vandoorn;resing;scheinhardt:98}[Section 1].\\
Because of technical difficulties originating from continuous state space, the authors discretize in
\cite{adan;vandoorn;resing;scheinhardt:98}[Section 4] the state of the buffer content: {\em Quanta of fluid} instead of {\em volume of fluid}. The size of the quantum is distributed exponentially.

\subsection{Finite state space}
\label{sect:OneSidedFinite}

Due to reversibility for classical birth-death processes the restriction from state space $\mathbb N_0$
to a finite state space $\{0, 1, \dots, N\}$ with $N\geq 1$ the stationary distribution of the latter one is obtained by truncation of the stationary distribution of the unrestricted process to the restricted
state  space.

This is not the case for the alternating birth-death process in a random environment, as can be seen by inserting 
 \eqref{eq:Oneuptheo1} and  \eqref{eq:Oneuptheo2} into the local balance equations.

For $n=0, 1,\dots,  N-1$ the equations \eqref{eq:GBOne1} through \eqref{eq:GBOne4}  fix again the relevant
conditions, and we have to add the boundary conditions for states $(N,b), (N,d)$. These are
\beqaro
x(N,b)\delta_N &=& x(N-1,b)\lambda_{N-1}+ x(N,d) \beta_N, \\ 
x(N,d)(\mu_N + \nu_N + \beta_N) &=&   x(N-1,b) \kappa_{N-1} + x(N,b) \delta_{N}. 
\eeqaro
Taking $x(N-1,b)$  from \eqref{eq:Oneuptheo1} these equations are  solved by
\beqo
x(N,b)= x(0,d)\beta_0 \big( \Lambda_{N-1} + \lambda_{N-1} \M_{N}\big)
\frac{\prod_{k = 1}^{N-1}\left(\Lambda_{k-1} \beta_{k} +\M_{k}\lambda_{k-1}\right)}
{\prod_{k = 0}^{N-1}\left(\Lambda_{k} \mu_{k+1}+\M_{k+1}\delta_{k}\right)}
\eeqo
and
\beqo
x(N,d) = x(0,d) \beta_0 \Lambda_{N-1}
\frac{\prod_{k = 1}^{N-1}\left(\Lambda_{k-1} \beta_{k} +\M_{k}\lambda_{k-1}\right)}
{\prod_{k = 0}^{N-1}\left(\Lambda_{k} \mu_{k+1}+\M_{k+1}\delta_{k}\right)}.
\eeqo
It is easy to see, that these solutions are compatible with \eqref{eq:GBOne1} - \eqref{eq:GBOne4}
and the solutions \eqref{eq:Oneuptheo1} and \eqref{eq:Oneuptheo2}.

\subsection{Regularity condition for the one-sided process}\label{sect:REGULAROneSided}

By definition, $Q$ it is conservative and under the assumption that the parameters
$\lambda_n, \kappa_n, \delta_n,$ $\mu_n, \nu_n, \beta_n$ are strictly positive it is irreducible.
Strict positivity is not necessary for irreducibility, but for simplicity of demonstration we shall put this assumption in force.
Because of the very general parameter set, $Q$ and an associated process in general are not regular, i.e.
under our assumptions a process constructed from $Q$ may explode in finite time with positive probability. The relevant criterion is ``Reuter's explosion condition'', see \cite{asmussen:03}[Proposition II.3.3].
\bprop\label{prop:reuter0} (Reuter (1957)) A Markovian jump process with discrete state space $E$ and transition intensity matrix $Q$ is nonexplosive if and only if  the set of  equations
\beq\label{eq:reuter0}
Q\cdot y = y
\eeq
has $y\equiv 0$ as the only nonnegative bounded solution. 
\eprop
Equation \eqref{eq:reuter0} is in our case with unknown
$y(n,b), y(n,d), n\geq 0$:
\beqar
y(0,b) (1+\lambda_0+\kappa_0+\delta_0) &=& y(0,d) \delta_0 + y(1,b)\lambda_0 + y(1,d)\kappa_0, \label{eq:reuter1}\\
y(0,d) (1+\beta_0) &=& y(0,b) \beta_0, \label{eq:reuter2}\\
y(n,b) (1+\lambda_n+\kappa_n+\delta_n) &=& y(n,d) \delta_n + y(n+1,b)\lambda_n + y(n+1,d)\kappa_n,
\quad n \geq 1, \label{eq:reuter3}\\
y(n,d) (1+\mu_n+\nu_n+\beta_n) &=& y(n,b) \beta_n + y(n-1,d)\mu_n + y(n-1,b)\nu_n.
\quad n \geq 1.\label{eq:reuter4}
\eeqar

\bprop\label{prop:ReuterOneSided}
An irreducible one-sided alternating birth-death process $\Z=(Z(t):t\geq 0)$ 
is non-explosive if and only if the partial solution sequence 
$(y(n,b):n\in \mathbb{N}_0)$ of the system \eqref{eq:reuter1} -  \eqref{eq:reuter4} increases unboundedly:\\
\beq\label{eq:div0}
y(n,b) \stackrel{n\to\infty}{\nearrow} \infty.
\eeq
Denote for $n\geq 0$ by $\Lambda^+_n := 1+\lambda_n + \kappa_n+\delta_n$  and
$\M^+_n := 1+\mu_n + \nu_n+\beta_n$ with $\mu_0 = \nu_0=0$.
For $n\geq 0$ holds
\beq\label{eq:Bounds0}
\prod_{k=0}^{n} 
\frac{(1+\lambda_k) \M^+_{k+1}+(1+\beta_{k+1})\kappa_k}
{\lambda_k \M^+_{k+1}+\beta_{k+1}\kappa_k} < y(n+1,b)
<\prod_{k=0}^{n} 
\frac{(1+\lambda_k+\delta_{k}) \M^+_{k+1}+(1+\mu_{k+1}+\beta_{k+1})\kappa_k}
{\lambda_k \M^+_{k+1}+\beta_{k+1}\kappa_k}
\eeq
A sufficient criterion for \eqref{eq:div0} is
\beq\label{eq:div0S}
\sum_{n=0}^{\infty} \frac{\M^+_{n+1}+\kappa_n }{\lambda_n \M^+_{n+1}+\beta_{n+1}\kappa_n} =\infty.
\eeq
\beq\label{eq:conv0}
\text{A sufficient criterion for\quad}\qquad\qquad
y(n,b) \stackrel{n\to\infty}{\nearrow} C <\infty\quad
\text{is}\quad\qquad\qquad\qquad
\eeq
\beq\label{eq:conv0S}
\sum_{n=0}^{\infty} \frac{(1+\delta_{n})\M^+_{n+1}+(1+\mu_{n+1})\kappa_n }{\lambda_n \M^+_{n+1}+\beta_{n+1}\kappa_n} <\infty.
\eeq
\eprop
\bproof
Because for irreducible $Q$ any nonnegative solution $y=(y(i):i\in E)$ of \eqref{eq:reuter0}, resp.
of \eqref{eq:reuter1} -  \eqref{eq:reuter4} with $y(i) =0$ for some $i\in E$ is identically zero, we can and will henceforth assume that $y(n,b)>0$ and $y(n,d)>0$ holds for all $n\in \mathbb{N}_0$.
We start with proving two facts for the solution of \eqref{eq:reuter0} for our process:
\begin{eqnarray}
&& (y(n,b):n\in\N_0)~~\text{is a strictly increasing sequence,}~~ 
	\text{and}\label{eq:increasing1}\\
&& y(n,b)-  y(n,d) > 0~~ \text{holds for all}~~ n\in\N_0. 
\label{eq:increasing2}
\end{eqnarray}
We combine \eqref{eq:reuter3} (resp. \eqref{eq:reuter1}) for $n\geq 0$ and \eqref{eq:reuter4} for the associated $n+1$, and obtain after some manipulations
\beqaro
y(n+1,b) =
y(n,b)\frac{\Lambda^+_n \M^+_{n+1}-\nu_{n+1}\kappa_n}{\lambda_n \M^+_{n+1}+\beta_{n+1}\kappa_n}
-y(n,d)\frac{\delta_n \M^+_{n+1}+\mu_{n+1}\kappa_n}{\lambda_n \M^+_{n+1}+\beta_{n+1}\kappa_n},
\\ 
y(n+1,d) =
y(n,b)\frac{(\Lambda^+_n \M^+_{n+1}-\nu_{n+1}\kappa_n)\beta_{n+1}+
	(\lambda_n \M^+_{n+1}+\beta_{n+1}\kappa_n)\nu_{n+1}}
{(\lambda_n \M^+_{n+1}+\beta_{n+1}\kappa_n)\M^+_{n+1}}\nonumber\\
-y(n,d)\frac{(\delta_n \M^+_{n+1}+\mu_{n+1}\kappa_n)\beta_{n+1}-
	(\lambda_n \M^+_{n+1}+\beta_{n+1}\kappa_n)\mu_{n+1}}
{(\lambda_n \M^+_{n+1}+\beta_{n+1}\kappa_n)\M^+_{n+1}}.
\eeqaro
Abbreviating
\beqo
D_n :=
\frac{\delta_n \M^+_{n+1}+\mu_{n+1}\kappa_n}{\lambda_n \M^+_{n+1}+\beta_{n+1}\kappa_n},\quad
E_n := \frac{(1+\lambda_n) \M^+_{n+1}+(1+\beta_{n+1})\kappa_n}
{\lambda_n \M^+_{n+1}+\beta_{n+1}\kappa_n},
\eeqo
and noticing 
\beqo
\frac{\Lambda^+_n \M^+_{n+1}-\nu_{n+1}\kappa_n}
{\lambda_n \M^+_{n+1}+\beta_{n+1}\kappa_n} = D_n + E_n,
\eeqo
this is for $n\geq 0$
\beqar
y(n+1,b) &=& y(n,b) (D_n+E_n) -y(n,d) D_n, \label{eq:reuter7}\\
y(n+1,d) &=& y(n,b) \left((D_n+E_n) \frac{\beta_{n+1}}{\M^+_{n+1}} + \frac{\nu_{n+1}}{\M^+_{n+1}}\right)
-y(n,d)\left( D_n \frac{\beta_{n+1}}{\M^+_{n+1}} - \frac{\mu_{n+1}}{\M^+_{n+1}}\right).\quad \label{eq:reuter8}
\eeqar
Note that $E_n>0, D_n>0$ holds. For later reference we remark that we obtained \eqref{eq:reuter7} and \eqref{eq:reuter8} without explicit reference to \eqref{eq:reuter2}.


We first remark that from \eqref{eq:reuter2} follows
\beq\label{eq:div2}
y(0,b) - y(0,d) > 0.
\eeq
Assume now that for some $n\geq 0$ holds
\beq\label{eq:div3}
y(n,b)-  y(n,d) > 0.
\eeq
From \eqref{eq:reuter7} and 
$E_n= 1 + \frac{\M^+_{n+1}+\kappa_n }{\lambda_n \M^+_{n+1}+\beta_{n+1}\kappa_n}$
it follows
\beq\label{eq:div4}
y(n+1,b) - y(n,b) =\underbrace{( y(n,b) - y(n,d) )}_{> 0,~\text{from }~\eqref{eq:div3} } D_n
+ \frac{\M^+_{n+1}+\kappa_n}{\lambda_n \M^+_{n+1}+\beta_{n+1}\kappa_n}y(n,b) > 0.
\eeq
Moreover, we have with \eqref{eq:reuter4}
\beqar
&&y(n+1,b)-  y(n+1,d)\label{eq:div4A}\\
&\stackrel{\eqref{eq:reuter4}}{=}& y(n+1,b) - \left(y(n+1,b) \frac{\beta_{n+1}}{\M^+_{n+1}} +
y(n,d)  \frac{\mu_{n+1}}{\M^+_{n+1}} + y(n,b)  \frac{\nu_{n+1}}{\M^+_{n+1}}\right)\nonumber\\
&=& y(n+1,b) \frac{1 + \mu_{n+1}+\nu_{n+1}}{\M^+_{n+1}} - y(n,b)  \frac{\nu_{n+1}}{\M^+_{n+1}}
- \underbrace{y(n,d)}_{< y(n,b) ~\text{from}~\eqref{eq:div3} }  \frac{\mu_{n+1}}{\M^+_{n+1}}\nonumber\\
&\geq&  \underbrace{\left((y(n+1,b) - y(n,b)\right)}_{> 0  ~\text{from}~\eqref{eq:div4} }
\frac{\mu_{n+1} + \nu_{n+1}}{\M^+_{n+1}} + y(n+1,b) \frac{1}{\M^+_{n+1}} >0.\nonumber
\eeqar
Summarizing, starting from \eqref{eq:div2}, we have proved 
that \eqref{eq:increasing1} (from \eqref{eq:div4}) and \eqref{eq:increasing2}  (from \eqref{eq:div4A}) hold.
Consequently, whenever $y=(y(n,t): n\in \N_0, t\in \{b,d\})$ is unbounded, then \eqref{eq:div0} holds (by \eqref{eq:increasing1}).
And clearly, vice versa.\\

To prove \eqref{eq:div0S} we rewrite \eqref{eq:reuter7} as
\beq\label{eq:div1S}
y(n+1,b)=y(n,b)\cdot E_n + D_n(y(n,b)-y(n,d)),
\eeq
which by iteration yields
\beq\label{eq:div2S}
y(n+1,b)= \Big(\underbrace{\prod_{k=0}^{n}E_k}_{(\star)}\Big)\cdot y(0,b)
+\sum_{m=0}^{n} \left(\prod_{k=m+1}^{n}E_k\right)\cdot D_m\cdot(y(m,b)-y(m,d)).
\eeq
From \eqref{eq:div2S} and \eqref{eq:increasing2} we obtain the left inequality of \eqref{eq:Bounds0}.
Moreover, the sequence of products $(\star)$ converges 
if and only if the series
\beqo
\sum_{n=0}^{\infty}  \frac{\M^+_{n+1}+\kappa_n }{\lambda_n \M^+_{n+1}+\beta_{n+1}\kappa_n} 
\eeqo
converges \cite{meschkowski:82}[Section IX.1, Kriterium XVII]. 
So, if that series diverges the sequence of products $(\star)$  diverges and (because of \eqref{eq:div4A}) the 
partial solution sequence $(y(n,b), n\in \mathbb{N}_0)$ 
diverges as well.\\
On the other side, from \eqref{eq:div1S} follows
\beq
y(n+1,b)<y(n,b)(E_n+D_n) \leq
\Big(\underbrace{\prod_{k=0}^{n}(E_k+D_k)}_{(\star\star)}\Big)\cdot y(0,b),
\eeq
which proves the right inequality of \eqref{eq:Bounds0}.
The sequence of products $(\star\star)$ converges 
if and only if the series
\beqo
\sum_{n=0}^{\infty}  \frac{(1+\delta_{n})\M^+_{n+1}+(1+\mu_{n+1})\kappa_n }{\lambda_n \M^+_{n+1}+\beta_{n+1}\kappa_n} 
\eeqo
converges by the same argument as for $(\star)$.
\eproof
\bremark
The iterative scheme \eqref{eq:reuter7} $+$ \eqref{eq:reuter8} provides a recursion to decide about
\eqref{eq:div0}:
\beqo
\begin{pmatrix}
	y(n+1,b)\\
	y(n+1,d)
\end{pmatrix}
=\begin{pmatrix}
D_n+E_n& -D_n\\
	(D_n+E_n)\cdot \frac{\beta_{n+1}}{\M^+_{n+1}}+ \frac{\nu_{n+1}}{\M^+_{n+1}}&
	-D_n  \cdot \frac{\beta_{n+1}}{\M^+_{n+1}}+ \frac{\nu_{n+1}}{\M^+_{n+1}}
\end{pmatrix}
\cdot \begin{pmatrix}
	y(n,b)\\
	y(n,d)
\end{pmatrix},
\eeqo
and we may choose without loss of generality as initial value
\beqo
\begin{pmatrix}
	y(0,b)\\
	y(0,d)
\end{pmatrix}
=
\begin{pmatrix}
	1\\
	{\beta_0}/{(1+\beta_0)}
\end{pmatrix},
\eeqo
\eremark
As can be seen by  careful inspection of the proof, strict positivity of the parameters is not necessary for the result to hold. Sufficient condition is irreducibility of $Q$.

\bremark\label{rem:NonMonotone1}
A consequence of \eqref{eq:div2} is that $(y(n,b):n\in\N_0)$ is a strictly increasing sequence.
A similar general property is not valid for $(y(n,d):n\in\N_0)$. This can be seen as follows.
Utilizing \eqref{eq:reuter2} and \eqref{eq:reuter8} we obtain
\beqaro
y(1,d)-y(0,d)
=y(0,d) \frac{1}{\M^+_{1}}\left(\frac{1}{\beta_0}\left\{(D_0+E_0) \beta_1+\nu_1 \right\} 
+\left[ \frac{\beta_1 \M^+_{1} + \beta_1\kappa_0}{\lambda_0 \M^+_{1} + \beta_1\kappa_0} - 1 \right]\right).
\eeqaro
Because $\left\{(D_0+E_0) \beta_1+\nu_1 \right\}>0$ this is 
 for $\lambda_0 \leq \beta_1$ strictly positive. On the other side, under $\lambda_0 > \beta_1$
the term in the squared brackets is in $(-1,0)$, and because $\left\{(D_0+E_0) \beta_1+\nu_1 \right\}>0$
is independent of $\beta_0$, by choosing $\beta_0$ appropriately, $y(1,d)-y(0,d)$ can be made positive or negative.
\eremark

\section{Two-sided alternating  birth-death process}
\label{sect:ModelTwoSided}
In this section we let the birth-death process take negative values as well, i.e. the reflecting boundary at $0$ is removed.
Let $\X=(X(t):t\geq 0)$ with $X(t): \omfp\to (\mathbb{Z},2^{\mathbb{Z}})$
be a variant of the standard two-sided birth-death processes the evolution of which depends on an external environment which changes randomly between two states $\{b, d\}$. We denote the environment process by
$\Y=(Y(t):t\geq 0)$ with $Y(t): \omfp\to (\{b, d\},2^{\{b, d\}})$.
Whenever the environment is in state $b$ (indicating that {\em births} may occur) the process $\X$ moves only upwards, by births occurring, and
whenever the environment is in state $d$ (indicating that {\em deaths} may occur) the process $\X$ moves only downwards, due to deaths occurring.

The joint process with state space $\mathbb{Z}\times \{b, d\}$ is denoted by $\Z=(\X,\Y)$ with
$Z(t) = (X(t),(Y(t))$. The movements of $\Z$ are governed by a transition intensity matrix
 $Q=(q(z,z'):z,z'\in \mathbb{Z}\times \{b, d\})$ with strictly positive entries as follows for $n\in \mathbb Z$:
\beqaro
q(n,b;n+1,b) &=& \lambda_n,~~\text{a birth occurs when the population is}~~n,\\ 
q(n,d;n-1,d) &=& \mu_n ,~~\text{a death occurs when the population is}~~n,\\ 
q(n,b;n,d) &=& \delta_n,~~\text{environment changes from ``birth'' to ``death'' when}\nonumber\\ &&\text{the population is}~~n, \\ 
q(n,d;n,b) &=& \beta_n,~~\text{environment changes from ``death'' to ``birth'' when}\nonumber\\
 &&\text{the population is}~~n,\\ 
q(n,b;n+1,d) &=& \kappa_n,~~\text{a birth occurs {\bf and} environment changes from ``birth''}\nonumber\\
 &&\text{ to ``death'' when the population is}~~n, \\ 
q(n,d;n-1,b) &=& \nu_n,~~\text{a death occurs {\bf and} environment changes from ``death''}\nonumber\\
 &&\text{to ``birth'' when the population is}~~n. 
\eeqaro
The diagonal elements of $Q$ are chosen such that row sums are zero. Unless otherwise indicated for special situations we assume throughout that the parameters
$\lambda_n, \kappa_n, \delta_n, \mu_n, \nu_n, \beta_n$ are strictly positive.\\
The transition graph for the two-sided process has the same structure as that of the one-sided version
as given in Figure \ref{fig:1-sided-transitions1}, now, for any $n\in \mathbb Z$.

\subsection{Stationary distribution for the two-sided process}\label{sect:GBTwoSided}
We assume in this section that the birth-death process is non-exploding, i.e. is regular. (In Section \ref{sect:REGULARTwoSided} we will discuss this in more detail.)
The global balance equations for $\Z$ are with unknown values $x(n,a)$ are for
$(n,a)\in \mathbb Z\times\{b,d\}$
\beqar
x(n,b)(\lambda_n + \kappa_n + \delta_n) &=& x(n-1,b)\lambda_{n-1} + x(n+1,d) \nu_{n+1}
+ x(n,d) \beta_n, \label{eq:GBTwo1}\\
x(n,d)(\mu_n + \nu_n + \beta_n) &=& x(n+1,d)\mu_{n+1} +  x(n-1,b) \kappa_{n-1} +
x(n,b) \delta_{n}. \label{eq:GBTwo2}
\eeqar

Again we observe that a sequence of simple cuts  reduces the computation to investigating a linear structure. We have for all  $n\in \mathbb Z$ the cut equation
\beqar
x(n,b)(\lambda_n + \kappa_n) &=&  x(n+1,d)(\mu_{n+1} + \nu_{n+1}), \label{eq:cutTwo1}\\
\text{which yields for all}~~ n\in \mathbb Z ~~~~~~~
x(n+1,d) &=& x(n,b)\frac{\Lambda_n}{\M_{n+1}} .  \label{eq:cutTwo2}
\eeqar
This allows to substitute in a first step the variables $x(n,d)$. In \eqref{eq:GBTwo1} we obtain
for $n\in \mathbb Z$
\beqo
x(n,b)(\lambda_n + \kappa_n + \delta_n) = x(n-1,b)\lambda_{n-1} +
[x(n,b)\frac{\Lambda_n}{\M_{n+1}}] \nu_{n+1}
+ [x(n-1,b)\frac{\Lambda_{n-1}}{\M_{n}}]  \beta_n,
\eeqo
and so
\beqo 
x(n,b)\frac{\Lambda_{n} \mu_{n+1}s+\M_{n+1}\delta_{n}}{\M_{n+1}}
= x(n-1,b)\frac{\lambda_{n-1}\M_n + \Lambda_{n-1} \beta_n}{\M_{n}},
\eeqo
which defines a two-sided standard birth-death process with state dependent birth and death rates.

For $n> 0$ this leads to
\beq
x(n,b)= x(n-1,b)\frac{\M_{n+1}}{\M_{n}}\cdot
\frac{\Lambda_{n-1} \beta_{n} +\M_{n}\lambda_{n-1}}{\Lambda_{n} \mu_{n+1}+\M_{n+1}\delta_{n}},
\label{eq:substTwo2}
\eeq
and iterating \eqref{eq:substTwo2} we obtain for all  $n> 0$
\beqaro
x(n,b)&=& x(0,b)\cdot \prod_{k = 1}^{n}\frac{\M_{k+1}}{\M_{k}}\cdot
\prod_{k = 1}^{n}
\frac{\Lambda_{k-1} \beta_{k} +\M_{k}\lambda_{k-1}}{\Lambda_{k} \mu_{k+1}+\M_{k+1}\delta_{k}}\nonumber\\
&=& x(0,b)\cdot \frac{\M_{n+1}}{\M_{1}}\cdot
\prod_{k = 1}^{n}
\frac{\Lambda_{k-1} \beta_{k} +\M_{k}\lambda_{k-1}}{\Lambda_{k} \mu_{k+1}+\M_{k+1}\delta_{k}}
\eeqaro

For $-n\leq 0$ we obtain
\beq
x(-n-1,b)= x(-n,b)\frac{\M_{-n}}{\M_{-n+1}}\cdot
\frac{\Lambda_{-n} \mu_{-n+1}+\M_{-n+1}\delta_{-n}}{\Lambda_{-n-1} \beta_{-n} +\M_{-n}\lambda_{-n-1}}.
\label{eq:substTwo3}
\eeq
Iterating \eqref{eq:substTwo3} we obtain for all  $-n< 0$
\beqaro
x(-n,b)&=& x(0,b)\cdot \prod_{k = -n}^{-1}\frac{\M_{k+1}}{\M_{k+2}}\cdot
\prod_{k = -n}^{-1}
\frac{\Lambda_{k+1} \mu_{k+2}+\M_{k+2}\delta_{k+1}}{\Lambda_{k} \beta_{k+1} +\M_{k+1}\lambda_{k}}\nonumber\\
&=& x(0,b)\cdot \frac{\M_{-n+1}}{\M_{1}}\cdot
\prod_{k = -n}^{-1}
\frac{\Lambda_{k+1} \mu_{k+2}+\M_{k+2}\delta_{k+1}}{\Lambda_{k} \beta_{k+1} +\M_{k+1}\lambda_{k}}.
\eeqaro

Taking into account \eqref{eq:cutTwo2} we obtain for all  $n> 0$
\beqo
x(n,d)= x(0,b)\cdot \frac{\Lambda_{n-1}}{\M_{1}}\cdot
\prod_{k = 1}^{n-1}
\frac{\Lambda_{k-1} \beta_{k} +\M_{k}\lambda_{k-1}}{\Lambda_{k} \mu_{k+1}+\M_{k+1}\delta_{k}},
\eeqo
and for $-n\leq 0$
\beqo
x(-n,d)= x(0,b)\cdot \frac{\Lambda_{-n-1}}{\M_{1}}\cdot
\prod_{k = -n-1}^{-1}
\frac{\Lambda_{k+1} \mu_{k+2}+\M_{k+2}\delta_{k+1}}{\Lambda_{k} \beta_{k+1} +\M_{k+1}\lambda_{k}}.
\eeqo
We summarize these computations as
\btheorem
Denote for $n\in\mathbb Z$ by $\Lambda_n := \lambda_n + \kappa_n$ the total uprate out of $n$ and $\M_n := \mu_n + \nu_n$ the total downrate out of $n$.
The global balance equation for $\Z$ are solved with any $x(0,b)>0$ by
\beqar
x(n,b)&=&  x(0,b)\cdot \frac{\M_{n+1}}{\M_{1}}\cdot
\prod_{k = 1}^{n}
\frac{\Lambda_{k-1} \beta_{k} +\M_{k}\lambda_{k-1}}{\Lambda_{k} \mu_{k+1}+\M_{k+1}\delta_{k}},
\quad n>0,
\label{eq:upTwo3}\\
x(-n,b)&=& x(0,b)\cdot \frac{\M_{-n+1}}{\M_{1}}\cdot
\prod_{k = -n}^{-1}
\frac{\Lambda_{k+1} \mu_{k+2}+\M_{k+2}\delta_{k+1}}{\Lambda_{k} \beta_{k+1} +\M_{k+1}\lambda_{k}}.
\quad n\geq 0,\label{eq:downTwo3}\\
x(n,d)&=& x(0,b)\cdot \frac{\Lambda_{n-1}}{\M_{1}}\cdot
\prod_{k = 1}^{n-1}
\frac{\Lambda_{k-1} \beta_{k} +\M_{k}\lambda_{k-1}}{\Lambda_{k} \mu_{k+1}+\M_{k+1}\delta_{k}},
\quad n>0,\label{eq:upTwo4}\\
x(-n,d)&=& x(0,b)\cdot \frac{\Lambda_{-n-1}}{\M_{1}}\cdot
\prod_{k = -n-1}^{-1}
\frac{\Lambda_{k+1} \mu_{k+2}+\M_{k+2}\delta_{k+1}}{\Lambda_{k} \beta_{k+1} +\M_{k+1}\lambda_{k}},
\quad n\geq 0.
\label{eq:downTwo4}
\eeqar

$\Z$ is ergodic if and only if
\beqaro
C(2)&=& \sum_{n=1}^\infty (\Lambda_n+\M_{n+1})\cdot
\prod_{k = 1}^{n}\frac{\left(\Lambda_{k-1} \beta_{k} +\M_{k}\lambda_{k-1}\right)}
{\left(\Lambda_{k} \mu_{k+1}+\M_{k+1}\delta_{k}\right)}\\
&&+ \sum_{n=1}^\infty (\Lambda_{-n}+\M_{-n+1})\cdot
\prod_{k = -n}^{-1}\frac{\left(\Lambda_{k+1} \mu_{k+2}+\M_{k+2}\delta_{k+1}\right)}
{\left(\Lambda_{k} \beta_{k+1} +\M_{k+1}\lambda_{k}\right)}
< \infty
\eeqaro
\etheorem

\subsubsection{Example: Telegraph process}\label{sect:Telegraph}

Kac investigated in 1974 a stochastic  model related to the telegrapher's equation \cite{kac:74}.
He considered a moving particle on the real line $\mathbb R$, starting at $0$, and being influenced by an alternating renewal process $\Y=(Y(t):t\geq 0)$ with exponential holding times.
The process $\Y$ with state space $\{\text{\sc l,r}\}$ with meaning $\{\text{\sc l:= left,r:= right}\}$ represents the environment for the moving particle 
and is Markov for its own:
Holding times for state {\sc l = left} are exponential($\beta$),
holding times for state {\sc r = right} are exponential($\delta$),
and all holding times are independent.\\
There is a one-way-interaction: 
The status of the environment determines the direction of the particle's movement.\\
$[\Y(\cdot) = \text{\sc l}]\Longrightarrow$ particle moves with constant velocity $v$ in direction $-\infty$\\
$[\Y(\cdot) = \text{\sc r}]\Longrightarrow$ particle moves with constant velocity $v$  in direction $+\infty$\\
Whenever $\Y$ jumps, the particle changes its direction immediately.\\
The particle's position at time $t\geq 0$ is denoted by $X(t)\in \mathbb R$.
Then $\Z=(\X,\Y)$ is assumed to be a Markov process.

Kac started his investigation by approximating the particle's position on a
lattice $\epsilon \cdot  \mathbb Z$ for $\epsilon\downarrow 0$ and in discrete time
with holding times  deterministic $=1$.
This generalized random walk in a random environment and its continuous limit
for $\epsilon\downarrow 0$ are called ``telegraph processes''.

The two-sided alternating birth-death process mimics the telegraph process in continuous time with exponential   holding times with the particle moving on $\mathbb Z$ when setting
\beqo
\forall n\in \mathbb Z:\quad \lambda_n =\eta,~~  \mu_n =\eta,~~ \beta_n=\beta,~~ \delta_n=\delta,~~
\kappa_n=0,~~ \nu_n=0.
\eeqo
This process is clearly not ergodic.\\
Interesting enough, the two-sided alternating birth-death process allows to construct generalized telegraph processes which are ergodic by controlling the speed of the particle
and by superposition of drifts. The mean speed of the particle in the general two-sided alternating birth-death process to move from $n$ to $n+1$ conditioned on $Y(\cdot)= b$ is $(\lambda_n+\kappa_n)$,
while the mean speed to move from $n$ to $n-1$ conditioned on $Y(\cdot)= d$ is $(\mu_n+\nu_n)$.
Telegraph processes with random speeds are investigated e.g. in \cite{stadje;zacks:04}, \cite{crimaldi;dicrescenzo;iuliano;martinucci:13}, \cite{degregorio:10}.

The natural generalized model of the telegraph process in terms of the two-sided alternating birth-death process is defined by the parameter setting
\beqo
\forall n\in \mathbb Z:\quad \lambda_n >0,~~  \mu_n >0,~~ \beta_n >0,~~ \delta_n >0,~~
\kappa_n=0,~~ \nu_n=0.
\eeqo
For the ergodic process the stationary distribution $\pi^T$ is with normalization constant $C(T)$
\beqar
\pi^T(n,{\text{\sc r}})& =&C(T)^{-1} \prod_{k = 1}^n\frac{\lambda_{k-1}}{\mu_{k}}
\prod_{k=1}^{n}\frac{\mu_{k}+\beta_{k}}{\lambda_{k}+\delta_{k}},
\quad n>0,\label{eq:TeleStat1}\\
\pi^T(n,{\text{\sc l}})& =&C(T)^{-1} \prod_{k = 1}^n\frac{\lambda_{k-1}}{\mu_{k}}
\prod_{k=1}^{n-1}\frac{\mu_{k}+\beta_{k}}{\lambda_{k}+\delta_{k}},
\quad n>0,\label{eq:TeleStat2}\\
\pi^T(-n,{\text{\sc r}})& =&C(T)^{-1} \prod_{k = -n}^{-1}\frac{\mu_{k+1}}{\lambda_{k}}
\prod_{k=-n}^{-1}\frac{\lambda_{k+1}+\delta_{k+1}}{\mu_{k+1}+\beta_{k+1}},
\quad n\geq 0,\label{eq:TeleStat3}\\
\pi^T(-n,{\text{\sc l}})& =&C(T)^{-1} \prod_{k = -n}^{-1}\frac{\mu_{k+1}}{\lambda_{k}}
\prod_{k=-n-1}^{-1}\frac{\lambda_{k+1}+\delta_{k+1}}{\mu_{k+1}+\beta_{k+1}},
\quad n\geq 0.\label{eq:TeleStat4}
\eeqar
The short notation
\beqo
\ell_{k}:=\frac{\lambda_{k}}{\lambda_{k+1}+\delta_{k+1}},\quad
m_{k}:=\frac{\mu_{k}}{\mu_{k}+\beta_{k}},\quad k\in\mathbb{Z},
\eeqo
reveals a birth-death process structure.
\beqar
\pi^T(n,{\text{\sc r}})& =&C(T)^{-1} 
\prod_{k = 1}^n\frac{\ell_{k-1}}{m_{k}}
\quad n>0,\label{eq:TeleStat1S}\\
\pi^T(n,{\text{\sc l}})& =&C(T)^{-1} 
\prod_{k = 1}^n\frac{\ell_{k-1}}{m_{k}}\cdot
\frac{\lambda_{n}+\delta_{n}}{\mu_{n}+\beta_{n}}
\quad n>0,\label{eq:TeleStat2S}\\
\pi^T(-n,{\text{\sc r}})& =&C(T)^{-1} 
\prod_{k = -n}^{-1}\frac{m_{k+1}}{\ell_{k}}
\quad n\geq 0,\label{eq:TeleStat3S}\\
\pi^T(-n,{\text{\sc l}})& =&C(T)^{-1} 
\prod_{k = -n}^{-1}\frac{m_{k+1}}{\ell_{k}}\cdot
\frac{\lambda_{-n}+\delta_{-n}}{\mu_{-n}+\beta_{-n}}
\quad n\geq 0.\label{eq:TeleStat4S}
\eeqar
The expressions  \eqref{eq:TeleStat2S} and \eqref{eq:TeleStat4S}
show that in case of a balanced system, i.e. the intensity to leave
state $(n,{\text{\sc r}})$ equals the intensity to leave
state $(n,{\text{\sc l}})$ for all $n\in\mathbb{Z}$, the probability for these states in equilibrium is the same.
Moreover,  the expressions  
\eqref{eq:TeleStat1}--\eqref{eq:TeleStat4S} show that in a balanced system the $\delta_{k}$ are linear functions of the 
$\beta_{k}$ (and vice versa), but they can vary over $(0,\infty)$ without changing the stationary distribution, which is now the stationary distribution of simple two-sided birth-death processes.\\

\textbf{Ergodic telegraph process with constant speed.}
As indicated above the classical telegraph process with constant speed, i.e. $\lambda_{k} = \mu_{k} =\eta > 0$  and with constant mean duration of the periods for traveling to the left, resp. right, i.e. $\beta_{k}=\beta, \delta_{k}=\delta$ (and $\kappa_k=\nu_k=0$) for all $k\in \mathbb{Z}$, is not ergodic.\\
We consider the situation of a telegraph process with constant speed  $\lambda_{k} = \mu_{k} =\eta > 0$ and given intensities 
$\beta_k, k\geq 0$ and $\delta_k, k\leq 0$ for the particle to finish an ongoing period of traveling towards zero and restarting its drift to $-\infty$, resp.  $\infty$.
Our aim is to control the process by selecting 
$\delta_k, k\geq 1,$
and $\beta_k, k\leq -1$, in such a way that the system stabilizes, i.e. the associated Markov process is ergodic.
For the ergodic process the stationary distribution $\pi^T$ would be with normalization constant  $C(T_\eta)$
\beqaro
\pi^T(n,{\text{\sc r}})& =&C(T_\eta)^{-1} 
\prod_{k=1}^{n}\frac{\eta+\beta_{k}}{\eta+\delta_{k}},\quad
\pi^T(n,{\text{\sc l}}) =C(T_\eta)^{-1} 
\prod_{k=1}^{n-1}\frac{\eta+\beta_{k}}{\eta+\delta_{k}},
\quad n>0, \\ 
\pi^T(-n,{\text{\sc r}})& =&C(T_\eta)^{-1} 
\prod_{k=-n}^{-1}\frac{\eta+\delta_{k+1}}{\eta+\beta_{k+1}},\quad
\pi^T(-n,{\text{\sc l}}) =C(T_\eta)^{-1} 
\prod_{k=-n-1}^{-1}\frac{\eta+\delta_{k+1}}{\eta+\beta_{k+1}},
\quad n\geq 0. 
\eeqaro
We have to construct sequences $\delta_k, k\geq 1,$ and 
$\beta_k, k\leq -1$, such that $C(T_\eta)<\infty$ holds, with
\beqaro
C(T_\eta)&=&
\underbrace{\sum_{n=1}^{\infty}
\prod_{k=1}^{n}\frac{\eta+\beta_{k}}{\eta+\delta_{k}}}_{(\star)}
+ \sum_{n=1}^{\infty}
\prod_{k=1}^{n-1}\frac{\eta+\beta_{k}}{\eta+\delta_{k}}\\
&+&1+\frac{\eta+\beta_{0}}{\eta+\delta_{0}}
+
\sum_{n=1}^{\infty}
\prod_{k=-n}^{-1}\frac{\eta+\delta_{k+1}}{\eta+\beta_{k+1}}
+ \underbrace{\sum_{n=1}^{\infty}
\prod_{k=-n-1}^{-1}\frac{\eta+\delta_{k+1}}{\eta+\beta_{k+1}}}_{(\star\star)}.
\eeqaro
The form of the summands in the infinite sums suggests to apply
a quotient test for the summands. Obviously it suffices to guarantee that the sums $(\star)$ and $(\star\star)$ are convergent. The recipe we shall apply is borrowed from the construction of the Bertrand-De Morgan test, see e.g. \cite{abramov:20}[Theorem 1]:\\
$\mathbf{(\star)}$ Take any sequence $r_n, n\geq 1,$ with $\liminf_{n\to\infty} r_n > 1$
and define on the positive axis the control intensities for restarting traveling to the left 
\beq\label{eq:ControlDelta1}
\delta_{n+1}:= \beta_{n+1} +(\frac{1}{n} + \frac{r_n}{n\ln n})\cdot
(\eta+\beta_{n+1}),\quad n\geq 1.
\eeq
$\mathbf{(\star\star)}$ Take any sequence $t_n, n\geq 1,$ with $\liminf_{n\to\infty} t_n > 1$ and define on the negative axis the control intensities for restarting traveling to the right
\beq\label{eq:ControlBeta1}
\beta_{-n-1}:= \delta_{-n-1} +(\frac{1}{n} + \frac{t_n}{n\ln n})\cdot
(\eta+\delta_{-n-1}),\quad  n\geq 1,
\eeq
A direct application of the Bertrand-De Morgan test \cite{abramov:20}[Theorem 1] guarantees convergence of both series. This test reads as follows.\\
If for  positive numbers $a_n, n\geq 1,$
we have for all sufficient large $n$ a representation
\[
\frac{a_n}{a_{n+1}} 
= 1 + \frac{1}{n} + \frac{s_n}{n\ln n},
~~\text{where}~~\liminf_{n\to\infty} s_n > 1,
\]
then the series $\sum_{n=1}^{\infty} a_n$ is convergent.\\
We apply the criterion to $\mathbf{(\star)}$ and $\mathbf{(\star\star)}$ with
\[
a_n:= \prod_{k=1}^{n}\frac{\eta+\beta_{k}}{\eta+\delta_{k}},
~\text{with}~ \delta_{k}~ \text{from}~ \eqref{eq:ControlDelta1},
~\text{resp.}~
a_n:= \prod_{k=-n-1}^{-1}\frac{\eta+\delta_{k}}{\eta+\beta_{k}}
~\text{with}~ \beta_{k}~ \text{from} ~\eqref{eq:ControlBeta1},
~~  n\geq 1.
\]
Note, that a refinement of the Bertrand-De Morgan test as given e.g. in  \cite{abramov:20}[Theorem 2] would produce refined control intensities in a similar way as we demonstrated here.

\subsection{Regularity condition for the two-sided process}\label{sect:REGULARTwoSided} 

Recall, that we assume the parameters
$\lambda_n, \kappa_n, \delta_n,$ $\mu_n, \nu_n, \beta_n$ to be strictly positive.
Because of the very general parameter set, $Q$ and an associated two-sided process in general may explode in finite time with positive probability. To find criteria for regularity we apply again Reuter's criterion from Proposition \ref{prop:reuter0}. Equation \eqref{eq:reuter0} now reads
\beqar
y(n,b) (1+\lambda_n+\kappa_n+\delta_n) &=& y(n,d) \delta_n + y(n+1,b)\lambda_n + y(n+1,d)\kappa_n, ~~ n\in \mathbb{Z}, \label{eq:reuter21}\\
y(n,d) (1+\mu_n+\nu_n+\beta_n) &=& y(n,b) \beta_n + y(n-1,d)\mu_n + y(n-1,b)\nu_n, ~~
n\in \mathbb{Z}. \label{eq:reuter22}
\eeqar
Because for irreducible $Q$ any nonnegative solution $y=(y(i):i\in E)$ of \eqref{eq:reuter0}, resp.
of \eqref{eq:reuter21} -  \eqref{eq:reuter22} with $y(i) =0$ for some $i\in E$ is identically zero, we can and will henceforth assume that $y(n,b)>0$ and $y(n,d)>0$ holds for all $n\in \mathbb{Z}$. We have to find conditions that guarantee that this solution is unbounded.\\ 
We extend the definitions  $\Lambda^+_n := 1+\lambda_n + \kappa_n+\delta_n$  and
$\M^+_n := 1+\mu_n + \nu_n+\beta_n$  now for $n\in \mathbb{Z}$.\\
Abbreviating for $n\geq 0$
\beqo
D_n :=
\frac{\delta_n \M^+_{n+1}+\mu_{n+1}\kappa_n}{\lambda_n \M^+_{n+1}+\beta_{n+1}\kappa_n},\quad
E_n := \frac{(1+\lambda_n) \M^+_{n+1}+(1+\beta_{n+1})\kappa_n}
{\lambda_n \M^+_{n+1}+\beta_{n+1}\kappa_n},
\eeqo
we obtain similar to the one-sided case in \eqref{eq:reuter7} and \eqref{eq:reuter8} for $n\geq 0$
\beqar
y(n+1,b) &=& y(n,b) (D_n+E_n) -y(n,d) D_n, \label{eq:reuter27}\\
y(n+1,d) &=& y(n,b) \left((D_n+E_n)\frac{\beta_{n+1}}{\M^+_{n+1}} + \frac{\nu_{n+1}}{\M^+_{n+1}}\right)
-y(n,d)\left( D_n \frac{\beta_{n+1}}{\M^+_{n+1}} - \frac{\mu_{n+1}}{\M^+_{n+1}}\right).~~~ \label{eq:reuter28}
\eeqar
Recall that $E_n>0, D_n>0$ holds.

\textbf{Defining an anchor.} Differently from the one-sided process
we have no information on the sign of $y(0,b)-y(0,d)$. If this would be positive, we could proceed as in Section \ref{sect:REGULAROneSided} to prove monotonicity properties of the set $\{y(n,b),y(n,d), n\in\N_0\}$.
In this case the vector $(y(0,b),y(0,d))$ serves as an anchor for computing the sequence of vectors $((y(n,b),y(n,d)), n\in\N_0)$.\\ 
From the structure of the system  
\eqref{eq:reuter21}-\eqref{eq:reuter22} it follows that if $y(k_0,b)-y(k_0,d)>0$  for some $k_0\in \mathbb{Z}$, we can elaborate on the sequence $((y(k,b),y(k,d)), k\geq k_0)$ analogously as in Section \ref{sect:REGULAROneSided} to obtain with anchor $(y(k_0,b),y(k_0,d))$ similar monotonicity results.\\
For convenience of readers we will henceforth assume that $(y(0,b),y(0,d))$ can serve as an anchor for our evaluation, but we will not prescribe $y(0,b)-y(0,d)>0$.
Transformation of the proofs, resp. results to different anchors is obvious in any case.

Starting from the anchor $(y(0,b),y(0,d))$, we obtain that for $-n\leq 0$ holds

\beqaro
y(-n-1,d) &=&
y(-n,d)\frac{\M^+_{-n}\Lambda^+_{-n-1} -\kappa_{-n-1}\nu_{-n}}
{\mu_{-n} \Lambda^+_{-n-1}+\delta_{-n-1}\nu_{-n}}
-y(n,b)\frac{\beta_{-n}\Lambda^+_{-n-1} +\lambda_{-n-1}\nu_{-n}}
{\mu_{-n} \Lambda^+_{-n-1}+\delta_{-n-1}\nu_{-n}}
\\ 
y(-n-1,b)& =&
y(-n,d)\left\{\frac{\M^+_{-n}\Lambda^+_{-n-1} -\kappa_{-n-1}\nu_{-n}}
{\mu_{-n} \Lambda^+_{-n-1}+\delta_{-n-1}\nu_{-n}}\cdot \frac{\delta_{-n-1}}{\Lambda^+_{-n-1}}
+\frac{\kappa_{-n-1}}{\Lambda^+_{-n-1}}\right\}
\\ 
&&-y(-n,b)\left\{\frac{\beta_{-n}\Lambda^+_{-n-1} +\lambda_{-n-1}\nu_{-n}}
{\mu_{-n} \Lambda^+_{-n-1}+\delta_{-n-1}\nu_{-n}} \cdot \frac{\delta_{-n-1}}{\Lambda^+_{-n-1}}
- \frac{\lambda_{-n-1}}{\Lambda^+_{-n-1}}\right\}
\nonumber
\eeqaro
Abbreviating
\beqo
C_{-n} := 
\frac{(1+\mu_{-n})\Lambda^+_{-n-1} +(1+\delta_{-n-1})\nu_{-n}}
{\mu_{-n} \Lambda^+_{-n-1}+\delta_{-n-1}\nu_{-n}}\quad
B_{-n} :=
\frac{\beta_{-n}\Lambda^+_{-n-1} +\lambda_{-n-1}\nu_{-n}}
{\mu_{-n} \Lambda^+_{-n-1}+\delta_{-n-1}\nu_{-n}}
\eeqo
and noticing 
\beqo
\frac{\M^+_{-n}\Lambda^+_{-n-1} -\kappa_{-n-1}\nu_{-n}}
{\mu_{-n} \Lambda^+_{-n-1}+\delta_{-n-1}\nu_{-n}}
=B_{-n} + C_{-n},
\eeqo
this is for $-n\leq 0$
\beqar
y(-n-1,d) &=&
y(-n,d)(B_{-n}+C_{-n}) -y(n,b)B_{-n}, \label{eq:reuter27Two}\\
y(-n-1,b) &=&
y(-n,d)\left\{(B_{-n}+C_{-n})\cdot \frac{\delta_{-n-1}}{\Lambda^+_{-n-1}}
+\frac{\kappa_{-n-1}}{\Lambda^+_{-n-1}}\right\}
\label{eq:reuter28}\\
&&\qquad\quad -y(n,b)\left\{B_{-n} \cdot \frac{\delta_{-n-1}}{\Lambda^+_{-n-1}}
- \frac{\lambda_{-n-1}}{\Lambda^+_{-n-1}}
\right\}.\nonumber
\eeqar
Then, starting from the anchor $(y(0,b),y(0,d))$ we obtain from  
\eqref{eq:reuter27} 
\beq\label{eq:div2STwo}
y(n+1,b)= \Big(\prod_{k=0}^{n}E_k\Big)\cdot y(0,b)
+\sum_{m=0}^{n} \left(\prod_{k=m+1}^{n}E_k\right)\cdot D_m\cdot(y(m,b)-y(m,d)),
\eeq
and from \eqref{eq:reuter27Two} symmetrically  via
$ y(-n-1,d) =
y(-n,d)C_{-n} +(y(-n,d)-y(n,b))B_{-n} 
:$
\beq\label{eq:div3STwo}
y(-n-1,d)= \Big(\prod_{k=0}^{n}C_{-k}\Big)\cdot y(0,d)
+\sum_{m=0}^{n} \left(\prod_{k=m+1}^{n}C_{-k}\right)\cdot B_{-m}\cdot(y(-m,d)-y(-m,b)),
\eeq
\blemma\label{lem:MonotoneTwo}
\textbf{(a)} If $y(0,b)-y(0,d) \geq 0,$ then
\begin{eqnarray}
&& (y(n,b):n\in\N_0)~~\text{is a strictly increasing sequence,}~~ 
\text{and}\label{eq:increasing1Two}\\
&& y(n,b)-  y(n,d) > 0~~ \text{holds for all}~~ n\geq 1. 
\label{eq:increasing2Two}
\end{eqnarray}
\textbf{(b)} If $y(0,d)-y(0,b) \geq 0,$ then
\begin{eqnarray}
&& (y(-n,d):n\in\N_0)~~\text{is a strictly increasing sequence,}~~ 
\text{and}\label{eq:increasing3Two}\\
&& y(-n,d)-  y(-n,b) > 0~~ \text{holds for all}~~ n\geq 1. 
\label{eq:increasing4Two}
\end{eqnarray}
\elemma
\bproof
\textbf{(a)} If $y(0,b)-y(0,d) > 0,$ the proof is verbatim the same as that for the parallel part of Proposition \ref{prop:ReuterOneSided}.\\
If $y(0,b)=y(0,d)>0,$ from \eqref{eq:reuter27} we have 
$y(1,b)=y(0,b) E_0$ and with 
$
E_n= 1  + \frac{\M^+_{n+1}+\kappa_n }
{\lambda_n \M^+_{n+1}+\beta_{n+1}\kappa_n}
$
follows $y(1,b) - y(0,b) >0$ and then analogously to \eqref{eq:div4A} we find $y(1,d) - y(1,b) >0$.
Restarting with the anchor $(y(1,b),y(1,d))$ we are back in the procedure of the one-sided case.\\
\textbf{(b)} Assume that $y(-n,d)-y(-n,b) > 0$ holds for some $n\geq 0$.
Taking into account 
\[
C_{-n} := 1 + \frac{\Lambda^+_{-n-1} +\nu_{-n}}
{\mu_{-n} \Lambda^+_{-n-1}+\delta_{-n-1}\nu_{-n}}
\]
we obtain from \eqref{eq:reuter27Two}
\beqo
y(-n-1,d) - y(-n,d) = (y(-n,d)-y(-n,b)) B_{-n} + 
y(-n,d)\frac{\Lambda^+_{-n-1} +\nu_{-n}}
{\mu_{-n} \Lambda^+_{-n-1}+\delta_{-n-1}\nu_{-n}} > 0,
\eeqo
and from \eqref{eq:reuter21}
\beqaro
&&y(-n-1,d) - y(-n-1,b)\\
& = &
(y(-n-1,d)-
\left(y(-n-1,d)\frac{\delta_{-n-1}}{\Lambda_{-n-1}^+} +
y(-n,b)\frac{\lambda_{-n-1}}{\Lambda_{-n-1}^+} +
y(-n,d)\frac{\kappa_{-n-1}}{\Lambda_{-n-1}^+}
\right)\\
&=&
y(-n-1,d)(1-\frac{\delta_{-n-1}}{\Lambda_{-n-1}^+}) -
\underbrace{y(-n,b)}_{< y(-n,d)}\frac{\lambda_{-n-1}}{\Lambda_{-n-1}^+} -
y(-n,d)\frac{\kappa_{-n-1}}{\Lambda_{-n-1}^+}\\
&\geq& y(-n-1,d)\frac{1+\lambda_{-n-1}+\kappa_{-n-1}}{\Lambda_{-n-1}^+} - y(-n,d)\frac{\lambda_{-n-1}}{\Lambda_{-n-1}^+}
- y(-n,d)\frac{\kappa_{-n-1}}{\Lambda_{-n-1}^+}\\
&=&(y(-n-1,d) - y(-n,d))\frac{\lambda_{-n-1}+\kappa_{-n-1}}{\Lambda_{-n-1}^+} +
y(-n-1,d)	\frac{1}{\Lambda_{-n-1}^+} > 0.
\eeqaro
If $y(-n,d)=y(-n,b)$ holds we obtain from \eqref{eq:reuter27Two}
\beqo
y(-1,d) - y(0,d) = y(0,d) 
\frac{\Lambda^+_{-1} +\nu_{0}}
{\mu_{0} \Lambda^+_{-1}+\delta_{-1}\nu_{0}} > 0,
\eeqo
and moreover
\beqo
y(-1,d) - y(-1,b) \geq (y(-1,d) -y(0,d))
\frac{\lambda_{-1} +\kappa_{-1}}{ \Lambda^+_{-1}} 
+y(-1,d)\frac{1}{ \Lambda^+_{-1}} > 0,
\eeqo
which yields a new anchor $(y(-1,d), y(-1,b))$ which satisfies
$y(-1,d) - y(-1,b) >0$ and we can restart the computations as in the first part of the proof of \textbf{(b)}.
\eproof

\bprop\label{prop:NonexplodingTwo} 
\textbf{(a)} A non-zero solution $((y(n,b),y(n,d)):n\in \mathbb{Z})$ of Reuter's regularity equation \eqref{eq:reuter0} is unbounded if and only if\\
\textbf{(i)} either the sequence $(y(n,b):n\in \mathbb{N}_0)$ is from some $n_+\geq 0$ on strictly increasing to $\infty$,\\
\textbf{(ii)} or the sequence $(y(-n,d):n\in \mathbb{N}_0)$ is from some $n_-\leq 0$ on strictly increasing to $\infty$,\\
\textbf{(iii)} or \textbf{(i)}  and \textbf{(ii)}  are valid concurrently.\\
\textbf{(b)} For the two-sided alternating birth-death process to be non-exploding (in finite time) it suffices that either of the series
\beq\label{eq:Non-ExplTwo1}
\sum_{n=0}^{\infty}  \frac{\M^+_{n+1}+\kappa_n }{\lambda_n \M^+_{n+1}+\beta_{n+1}\kappa_n} 
\eeq
or the series
\beq\label{eq:Non-ExplTwo2}
\sum_{n=0}^{\infty}  \frac{\Lambda^+_{-n-1}+\nu_n} {\mu_{-n}\Lambda^+_{-n-1}+\delta_{-n-1}\nu_n} 
\eeq
is divergent (or both).
\eprop
\bproof
We start with proving part \textbf{(b)}.
We can find an anchor for our computations as follows. 
Consider a (strict) positive solution $y:=((y(n,b),y(n,d)): n\in\mathbb{Z})$ of Reuter's regularity equation \eqref{eq:reuter0}
and denote
\beq\label{eq:Anchor1}
A := \inf\{k\in\mathbb{Z}: y(k,b) - y(k,d)\geq 0\}.
\eeq
\textbf{Case 1.} $A= -\infty.$ Then for all $n\in\mathbb{Z})$ holds
$y(n,b) - y(n,d) >0$. Taking $(y(0,b),y(0,d))$ as anchor we conclude with Lemma \ref{lem:MonotoneTwo}\textbf{(a)} that the sequence $(y(n,b):n\in\mathbb{N}_0)$ is strictly increasing and all differences $y(n,b) - y(n,d), n\in\mathbb{N}_0,$ are strictly positive.
From \eqref{eq:div2STwo} we see that divergence of the sequence of products $\prod_{k=0}^{n}E_{k}$ guarantees that the solution 
$y:=((y(n,b),y(n,d)): n\in\mathbb{Z})$ of equation \eqref{eq:reuter0} is unbounded.
This sequence of products converges if and only if
the series
\beqo
\sum_{n=0}^{\infty}  \frac{\M^+_{n+1}+\kappa_n }{\lambda_n \M^+_{n+1}+\beta_{n+1}\kappa_n} 
\eeqo
converges \cite{meschkowski:82}[Section IX.1, Kriterium XVII]. 
So, if that series diverges the sequence of products diverges and  the partial solution sequence $(y(n,b), n\in \mathbb{N}_0)$ 
diverges as well.\\
\textbf{Case 2.} $A= \infty.$ Then for all $n\in\mathbb{Z}$ holds
$y(n,d) - y(n,b) >0$. Taking $(y(0,b),y(0,d))$ as anchor we conclude with Lemma \ref{lem:MonotoneTwo}\textbf{(b)} that the sequence $(y(-n,d):n\in\mathbb{N}_0)$ is strictly increasing and all differences $y(-n,d) - y(-n,b), n\in\mathbb{N}_0,$ are strictly positive.
From \eqref{eq:div3STwo} we see that divergence of the sequence of products $\prod_{k=0}^{n}C_{-k}$ guarantees that the solution 
$y:=((y(n,b),y(n,d)): n\in\mathbb{Z})$ of equation \eqref{eq:reuter0} is unbounded.
This sequence of products converges if and only if
the series
\beq  \label{eq:Non-ExplTwo2Case2}
\sum_{n=0}^{\infty}  \frac{\Lambda^+_{-n-1}+\nu_n} {\mu_{-n}\Lambda^+_{-n-1}+\delta_{-n-1}\nu_n} 
\eeq
converges \cite{meschkowski:82}[Section IX.1, Kriterium XVII]. 
So, if that series diverges the sequence of products diverges and  the partial solution sequence $(y(-n,d), n\in \mathbb{N}_0)$ 
diverges as well.\\
\textbf{Case 3.} $A\in\mathbb{Z}$ is finite.
Then for all $n< A$ holds
$y(n,d) - y(n,b) >0$. Taking $(y(A-1,b),y(A-1,d))$ as anchor we conclude similar to Lemma \ref{lem:MonotoneTwo}\textbf{(b)} that the sequence $(y(n,d):n\leq A-1)$ is strictly increasing and for all $n\leq A-1$ the differences $y(n,d) - y(n,b)$  are strictly positive.
If the series
\beq \label{eq:Non-ExplTwo2Case3}
\sum_{n=0}^{\infty}  \frac{\Lambda^+_{-n-1}+\nu_n} {\mu_{-n}\Lambda^+_{-n-1}+\delta_{-n-1}\nu_n} 
\eeq
diverges we conclude as in Case 2. that the solution of
\eqref{eq:reuter0} is unbounded.\\
If \eqref{eq:Non-ExplTwo2Case3} is convergent and 
\eqref{eq:Non-ExplTwo2Case2} is divergent, we start with
$(y(A,b),y(A,d))$ as anchor and conclude similar to Lemma \ref{lem:MonotoneTwo}\textbf{(a)} that the sequence $(y(n,b):n\geq A)$ is strictly increasing and all differences $y(n,b) - y(n,d), n\geq A,$ are strictly positive. Arguing as in Case 1. we finish the proof.\\
For part \textbf{(a)} we notice that a non-zero solution of
\eqref{eq:reuter0}
is unbounded if and only if at least one of the following conditions is fulfilled.
\beqo
\limsup_{n \geq 0} y(n,b) = \infty,\quad
\limsup_{n\geq 0} y(n,d) = \infty,\quad
\limsup_{n\geq 0} y(-n,b) = \infty,\quad  
\limsup_{n\geq 0} y(-n,d) = \infty.
\eeqo 
Combining the facts from Lemma \ref{lem:MonotoneTwo} and from the proof of part \textbf{(b)} we conclude the statement of \textbf{(a)}.  
\eproof
\bremark
\textbf{(a)} The proof of Proposition \ref{prop:NonexplodingTwo} shows implicitly that the value $A$ in \eqref{eq:Anchor1} is uniquely defined.\\
\textbf{(b)} Evaluating bounds for the  $y(-n,b)$,
resp. $y(-n,d)$, similar to the one-sided case in Proposition 
\ref{prop:ReuterOneSided} is along the same lines of computation as
presented in Section \ref{sect:REGULAROneSided}.
\eremark

{\bf Acknowledgement:} I thank Jacques Resing for introducing to me his work in
\cite{adan;resing:96} and \cite{adan;vandoorn;resing;scheinhardt:98}, and
Vyacheslav Abramov for sending me his forthcoming paper \cite{abramov:20}.

\bibliographystyle{alpha}

\end{document}